    \documentclass[11pt]{extarticle}
        \usepackage[mathscr]{euscript}
        \usepackage[utf8]{inputenc}
        \usepackage[lowercase]{theoremref}
        \usepackage{ amsmath, amsfonts, amsthm, amssymb, csquotes, centernot, enumerate, tools, calc, soul}
        \makeatletter
            \newcommand{\mymargins@@@}{1in}
            \usepackage[left=\mymargins@@@, right=\mymargins@@@, top=\mymargins@@@, bottom=\mymargins@@@]{geometry}
        \makeatother
        \makeatletter
            \newcommand{\myTitle}{?}
            \newcommand{\headerSize@}{\footnotesize}
            \usepackage{fancyhdr}     
            \newcommand{\Begin}[1][] { \renewcommand{\myTitle}{#1} 
\begin{document} }
            \newcommand{\End}{\end{document}}
        \makeatother
    %}
    \usepackage{cite}
    \bibliographystyle{plain}
    
    \numberwithin{equation}{section}
    \linespread{1}
    \newcommand{\verbose}[1]{}
    
%}

%------------------------------------------------------------------

        %NOTE: There are still results from the interpolation section which are used in the sampling section. I am considering factoring these to their own section.
        
        %NOTE: I have tried using the following conventions:
            %(1) $w, z, \zeta$ for points in $\dd$;
            %(2) $p,q,s$ for function exponents;
            %(3) $\alpha, \lambda, \lambda$ for weight exponents;
            %(4) $b,c$ for real parameters;
            %(5) $S$ for operators;
            %(6) $\mu$ for measures.
            
%------------------------------------------------------------------

\providecommand{\keywords}[1]
{
  \small	
  \textbf{\textit{Keywords:}} #1
}

\title{Interpolation and Sampling in Analytic Tent Spaces}

\author{Caleb Parks}

\Begin[Interpolation and Sampling in Analytic Tent Spaces]

\maketitle\thispagestyle{empty}
\begin{abstract}
Introduced by Coifman, Meyer, and Stein, the tent spaces have seen wide applications in Harmonic Analysis. Their analytic cousins have seen some applications involving the derivatives of Hardy space functions. Moreover, the tent spaces have been a recent focus of research. We introduce the concept of interpolating and sampling sequences for analytic tent spaces analogously to the same concepts for Bergman spaces. We then characterize such sequences in terms of Seip's upper and lower uniform density. We accomplish this by exploiting a kind of M\"obius invariance for the tent spaces.

\end{abstract}


\section{Introduction}

Given a space $X$ of analytic functions on the unit disc, $\dd$, and a sequence $\set{z_k}_k$ in $\dd$, we consider the sequence-valued evaluation operator on $X$ defined by
$$
    Ef = \set{f(z_k)}_k
$$
From here, we can try to determine a sequence space $\mathscr{X}$ such that $E:X \arrow \mathscr{X}$ is a bounded operator. The space $\mathscr{X}$ is generally derived from $X$ in some fashion. In the case of the Hardy space, $H^p$, the space $\mathscr{X}$ is the weighted sequence space 
$$
    \ell_{-1}^{p}(\set{z_k}) = \set*{ \set{a_k}_k : \sum \abs{a_k}^p(1-\abs{z_k}^2) < \infty } \,.
$$
Please refer to page 149 of \cite{Duren} for more exposition. In the case of the Bergman space,
$$
    A_{\alpha}^{p} = \set*{ f \in \holo(\dd) : C_{\alpha} \int_{\dd} \abs{f(z)}^p (1-\abs{z}^2)^{\alpha} \dif A(z) < \infty }\,,
$$
where $\holo(\dd)$ is the set of holomorphic functions on $\dd$, the relevant space $\mathscr{X}$ is the sequence space
$$
    \ell_{\alpha}^{p}(\set{z_k}) = \set*{ \set{a_k}_k : \sum \abs{a_k}^p(1-\abs{z_k}^2)^{2+\alpha} < \infty } \,.
$$
For this result, see section 7.2 of \citePage{166}{Zhu}. For more discussion on how to choose such $\mathscr{X}$, see \citePage{15}{Seip_Book}, \citePage{41}{Seip_Book}, \citePage{64}{Seip_Book}, and \citePage{716}{Seip_Survey}. Once $\mathscr{X}$ has been fixed, two properties of $E$ are given special names: 
\begin{itemize}
    \item If $E:X \arrow \mathscr{X}$ is surjective, we say $\set{z_k}$ is \emph{interpolating} for $(X,\mathscr{X})$.
    \item If $E:X \arrow \mathscr{X}$ is injective with closed range, we say $\set{z_k}$ is \emph{sampling} for $(X,\mathscr{X})$. 
\end{itemize}
When there is a natural choice for $\mathscr{X}$, we drop $\mathscr{X}$ from the notation. Many cases of these questions have been studied. The cases where $X$ is the Hardy space  (\cite{Carleson}, \cite{ShapiroShields}, \cite{Kabaila}), the Bergman Space (\cite{Seip_InterpolationSampling}, \cite{Schuster_Thesis}, \cite{Schuster_SamplingBig}, \cite{Schuster_SamplingSmall}, \cite{Schuster_Surfaces}), the Bloch space \cite{BoeNicolau}, the Fock Space \cite{SeipOrtega_Fock}, and the Paley-Weiner space \citePage{98}{Seip_Book} have all been studied. 

In the cases studied, people have connected interpolation to many results. In the Hardy space, interpolation can be used in proving the \enquote{corona theorem}; see \citePage{203}{Duren} and \cite{Carleson}. Also, the invertibility of certain Toeplitz operators $T:H^2 \to H^2$ is related to interpolation and sampling for a related space; see \citePage{88}{Seip_Book}. In Bergman spaces, Hedenmalm, Richter and Seip found applications for interpolation to construct invariant subspaces of infinite index; see \cite{HedenmalmRichterSeip}. Also in Bergman spaces, sampling sequences give a more precise version of the classical atomic decomposition \citePage{195}{DurenSchuster}. Finally, in \cite{Thomson} Thomson applies Seip's sampling theorem to show that the closure of the polynomials in $L^p(\dd,\dif\mu)$ can change quite dramatically with the parameter $p$; see \citePage{172}{HedenmalmKorenblumZhu} for more.

As we have noted the utility of sampling and interpolating sets, we hope to extend these results to the case of the analytic tent spaces. Let $\mu$ be a Borel measure and $L^0$ be the set of Borel-measurable functions from $\dd$ to $\cc$. Introduced in \cite{CoifmanMeyerStein} for the upper-half space, the \emph{tent spaces} are defined by
$$
    T_{\alpha}^{p,q}(\mu) = \set*{ f \in L^0 \myspace[1] \mid \myspace[1] \norm{f}_{T_{\alpha}^{p,q}(\mu)} < \infty }
$$
where the norm is given by
$$
    \norm{f}_{T_{\alpha}^{p,q}(\mu)}^p \defeq \int_{\bd\dd} \group{ \int_{\Gamma_{\zeta}} \abs{f(z)}^q \dif\mu_{\alpha-1} }^{p/q} \dif l(\zeta)
$$
with $\dif \mu_{\alpha} = (1-\abs{z}^2)^{\alpha} \dif\mu(z)$, $l$ being normalized arc-length measure, and $\Gamma_{\zeta}$ being a nontangential approach region defined by
$$
    \Gamma_{\zeta} = \conv\group[\big]{ \set[\big]{ z \in \cc : \abs{z} < 1/\sqrt{2} } \cup \zeta } \cap \dd\, .
$$
where $\conv(E)$ denotes the convex hull of $E$. The \emph{analytic tent space}, $AT_{\alpha}^{p,q}(\mu),$ is the intersection of $T_{\alpha}^{p,q}(\mu)$ with the space of analytic functions.
In the case that $\mu$ is normalized area measure, $A$, we write $T_{\alpha}^{p,q}$ and $AT_{\alpha}^{p,q}$ instead of $T_{\alpha}^{p,q}(A)$ and $AT_{\alpha}^{p,q}(A)$. We also define $\sigma = \dif A_{-1}$. In this notation, 
$$
\norm{f}_{T_{\alpha}^{p,q}}^p \defeq \int_{\bd\dd} \group{ \int_{\Gamma_{\zeta}} \abs{f}^q \dif\sigma_{\alpha} }^{p/q} \dif l(\zeta) \,.
$$
The choice of normalization for $T_{\alpha}^{p,q}$ gives that $T_{\alpha}^{p,p} = L_{\alpha}^p$ where $L_{\alpha}^p = L^p(A_{\alpha})$ is the Lebesgue space. Sometimes, it will help to define
$$
    \scriptA_{\alpha}^q[f;\mu](\zeta) = \group{\int_{\Gamma_{\zeta}} \abs{f}^q \dif \mu_{\alpha-1} }^{1/q} \,.
$$
When $\mu = A$, we will drop the dependence on $\mu$. 

The tent spaces have been a tool to study Hardy spaces. The classical Hardy spaces $H^p$ have derivatives which are in the tent space $T_1^{p,2}$. The tent spaces have seen applcations in the following papers: Luecking, \cite{Luecking_Embedding}, formulated Carleson-style embedding theorems, Zhao, \cite{Zhao}, characterized multiplication operators from the Hardy space $H^p$ to the weighted Bergman space $A_{\alpha}^q$ (see \cite{Zhao}), Miihkinen, Pau, Per\"al\"a, and Wang, \cite{MiihkinenPauPerala}, characterized the Volterra operators from Bergman spaces to Hardy spaces, and finally Cascante and Ortega, \cite{CascanteOrtega}, were able to characterize pointwise multipliers for the Hardy-Sobolev space.

The tent spacess have also been a recent focus of independent research. In the papers \cite{Amenta_Homogeneous, Amenta_Interpolation}, Alex Amenta extended the definition of tent spaces and developed the theory of weighted spaces in the non-holomorphic setting. In \cite{CohnVerbitsky}, Cohn and Verbitsky obtained a general factorization method for tent spaces; they applied this factorization to Hankel operators. Pel\'aez, R\"atty\"a, and Sierra have also studied embeddings of Bergman spaces into tent spaces (see \cite{PelaezRattyaSierra}). 

In harmonic analysis and differential equations, the tent spaces have also many applications. In their recent paper \cite{BartonHofmannMayboroda}; Barton, Hofmann, and Mayboroda studied solutions of certain differential operators with bounds on the $T_{\alpha}^{p,2}$-norm. Aucher and Prisuelos-Arribas have recently proved boundedness of certain maximal operators, Calder\'on-Zygmund operators, and Riesz potentials in tent spaces (see \cite{AuscherArribas}). In an attempt to study more general elliptic differential operators $L$, Hofmann and Mayboroda used the tent spaces to define a version of the real hardy space $H_L^p$ that is \enquote{adapted} to $L$. See the following for applications of theses spaces: \cite{HofmannMayboroda}, \cite{HofmannMayborodaMcintosh}\cite{KunstmannUhl}. 

For a sequence $\set{ z_k } \subseteq \dd$, we define $t_{\alpha}^{p,q} = T_{\alpha}^{p,q}(\nu^{\set{z_k}})$ where $\nu^{\set{z_k}}$ is defined by 
$$
    \dif \nu^{\set{z_k}} = \sum_k (1-\abs{z_k}^2)^2 \dif\delta_{z_k}
$$
and $\delta_{z_k}$ is the point-mass measure at $z_k$. The elements $f \in t_{\alpha}^{p,q}$ can be thought of as sequences with $a_k = f(z_k)$. Hence, we will identify $t_{\alpha}^{p,q}$ with the set of sequences $a_k$ such that
$$
    \norm{\set{a_k}}_{t_{\alpha}^{p,q}(\set{z_k})}^p \defeq 
    \int_{\bd\dd} \group{ \sum_{z_k \in \Gamma_{\zeta}} \abs{a_k}^q(1-\abs{z_k}^2)^{\alpha +1} }^{p/q}
$$
is finite. We will call $t_{\alpha}^{p,q}$ a \emph{sequence space of tent type}. As $t_{\alpha}^{p,q}$ is the tent-space analogue of $\ell_{\alpha}^p$, we feel it is appropriate to take $E:AT_{\alpha}^{p,q} \to t_{\alpha}^{p,q}$. After our discussion of interpolation sequences, we will extend the notion of interpolation for measures. We will prove that this extended notion does not yield anything new in the case of Carleson measures. We do not know what happens for non-Carleson measures.

As in the Bergman space, the \emph{pseudohyperbolic metric} $\rho$ is the correct distance to use when working with the tent spaces. We define 
$$
    \rho(z,w) = \abs{ \phi_w(z) }
$$
for each $z,w \in \dd$ where 
$$
    \phi_w(z) = \frac{ w-z }{ 1-\overbar{w}z }
$$
defines the unique involutory automorphism of the unit disc which interchanges $w$ and $0$. We denote the pseudohyperbolic disc centered at $z$ of radius $r$ by $D(z,r)$. We say that a sequence $\set{z_k}$ is \emph{uniformly discrete} if there is $\delta>0$ such that $\rho(z_j,z_k) \geq \delta$ for all $j \neq k$; we say that the largest such $\delta$ is the \emph{separation constant} denoting it by $\delta(\set{z_k})$. 

We now define the densities with which we can state the main theorem. Let $\set{z_k}$ be a sequence. For $r>1/2$ and $w\in \dd$, we define the \emph{predensity} $D(\set{z_k},w,r)$ by the equation
$$
    D(\set{z_k},w,r) = \bracket[\big]{\log(1-r)}^{-1} \sum_{ z_k \in \Hat{D}(w,r)} \log \abs{\phi_w(z_k) }
$$
where $\Hat{D}(w,r) = D(w,r) \excise D(w,1/2)$. We now define the \emph{upper uniform density} of $\set{z_k}$ by
$$
    D^+(\set{z_k}) = \limsup_{r \arrow 1} \group{ \sup_{w \in \dd} D(\set{z_k},w,r) }
$$
and the \emph{lower uniform density} of $\set{z_k}$ by
$$
    D^-(\set{z_k}) = \liminf_{r \arrow 1} \group{ \inf_{w \in \dd} D(\set{z_k},w,r) } \,.
$$ 

To be more concrete, we give explicit definitions of interpolation and sampling for our case: Let $\set{z_k}$ be a sequence in $\dd$. We say that $\set{z_k}$ is interpolating for $AT_{\alpha}^{p,q}$ if for each $\set{a_k} \in t_{\alpha}^{p,q}(\set{z_k})$ there is some $f \in AT_{\alpha}^{p,q}$ such that $a_k = f(z_k)$ for all $k$. We say that $\set{z_k}$ is sampling for $AT_{\alpha}^{p,q}$ if the inequalities
$$
    (1/C) \norm{f}_{T_{\alpha}^{p,q}} \leq \norm{f}_{t_{\alpha}^{p,q}(\set{z_k})} \leq C \norm{f}_{T_{\alpha}^{p,q}}
$$
hold for $C$ independent of $f$. Our main theorems are as follows:
\MainTheorem{interpolationTheorem}{
    The sequence $\set{z_k} \subseteq \dd$ is interpolating for $AT_{\alpha}^{p,q}$ if and only if $z_k$ is uniformly discrete and $D^+(\set{z_k}) < (1+\alpha)/q$.
}
\MainTheorem{samplingTheorem}{
    The sequence $\set{z_k} \subseteq \dd$ is sampling for $AT_{\alpha}^{p,q}$ if and only if $z_k$ is a finite union of uniformly discrete sequences and there is a uniformly discrete subsequence $\set{z_k'}$ of $\set{z_k}$ such that $D^-(\set{z_k'}) > (1+\alpha)/q$.
}
\setcounter{mainTheorem}{0}
For sake of comparison, we now recall the theorems of Seip for the Bergman spaces: Let $\set{z_k}$ be a sequence in $\dd$. We say that $\set{z_k}$ is interpolating for $A_{\alpha}^{p}$ if for each $\set{a_k} \in \ell_{\alpha}^{p}(\set{z_k})$ there is some $f \in A_{\alpha}^{p}$ such that $a_k = f(z_k)$ for all $k$. We say that $\set{z_k}$ is sampling for $A_{\alpha}^{p}$ if the inequalities
$$
    (1/C) \norm{f}_{A_{\alpha}^{p}} \leq \norm{f}_{\ell_{\alpha}^{p}(\set{z_k})} \leq C \norm{f}_{A_{\alpha}^{p}}
$$
hold for $C$ independent of $f$. Seip's theorems are as follows:
\SpecialTheorem{?A?}{
    A sequence $\set{z_k} \subseteq \dd$ is interpolating for the Bergman space $A_{\alpha}^p$ if and only if $\set{z_k}$ is uniformly discrete and $D^+(\set{z_k}) < (1+\alpha)/p$.
}

\SpecialTheorem{?B?}{
    A sequence $\set{z_k} \subseteq \dd$ is sampling for the Bergman space $A_{\alpha}^p$ if and only if $\set{z_k}$ is a finite union of uniformly discrete sequences and there is a uniformly discrete subsequence $\set{z_k'}$ of $\set{z_k}$ with $D^-({z_k'}) > (1+\alpha)/p$.
}
Please note the similarities between Seip's theorems and our own. In fact, we will use \thref{?A?} to prove one implication of \thref{interpolationTheorem}. 

\section{Definitions, Notations, and Background}
We now make a few assumptions for the duration of the paper unless otherwise specified: $p$ and $q$ denote positive, finite real numbers. $\alpha$ denotes a real number greater than $-1$. Any constants are always allowed to depend on $p$, $q$, and $\alpha$ unless otherwise specified. All sequences will be assumed infinite unless otherwise specified. 

This paper is primarily motivated to determine the behavior of the values of function $f$ at a sequence $\set{z_k}_k$ given that $f$ is in a given norm space. At some point, we will need to consider sequences of sequences. As such, we will think of a sequence as a function $\nn \arrow \cc$ and use the symbol $\Dott$ to denote its argument. Thus we will identify the function $z\0$ with the sequence $\set{z_k}_k$.  Moreover, given a sequence $z\0$ and a continuous function $f$, we will define the sequence $f(z\0)$ by $f(z\0)_k = f(z_k)$ whenever $z\0$ is in the domain of $f$. Thus we extend $f$ to map sequences to sequences in the natural way.

We will think of a tuple on an index set $J$ with values in $S$ as being a function $E:J \arrow S$. We will write $E_x$ for $E(x)$. A \emph{rearrangement} of a tuple $E:J \arrow S$ is the tuple $E \comp \phi$ for some bijection $\phi:J \arrow J$. We say two tuples $E$ and $E'$ are \emph{rearrangements} of each other if $E = E' \comp \phi$ for some rearrangement $\phi$. We define a \emph{set with repetition} as a tuple modulo rearrangement. Given two sets with repetition, $E:J \arrow S$ and $E':J' \arrow S$, we define
$$
    F: J \amalg J' \arrow S
$$
by
$$
F(x) =
\begin{cases}
    E(x) &: x\in J \\
    E'(x) &: x\in J'
\end{cases}
\,,
$$
modulo rearrangement and write $F = E \cup E'$. We define the union of sequences as an enumeration of $z\0 \cup w\0$ where $z\0$ and $w\0$ are thought of as sets with repetition. We also write $z\0 \subseteq E$ to mean that the image of $z\0$ lies inside of $E$ and $z \in z\0$ to mean $z$ is in the image of $z\0$. 

We now discuss the definition of the tent spaces in more detail. We also discuss some basic properties of tent spaces. We first have the following definitions. Given $r \in (0,1)$, we define a \emph{Stoltz family} of \emph{aperture} $r$ to be a family of sets 
$$
    \Lambda_{\zeta}^r = \conv( D(0,r) \cup \set{\zeta})
$$
for $\zeta \in \bd\dd$. Each element $\Lambda_{\zeta}^r$ of the family is called a \emph{Stoltz region} or \emph{Stoltz angle}. If we do not need to specify the aperture, we write $\Lambda_{\zeta}$ instead of $\Lambda_{\zeta}^r$. Distinguishing the case where $r = 1/\sqrt{2}$, we define $\Gamma_{\zeta} = \Lambda_{\zeta}^{1/\sqrt{2}}$. 

We will see soon that 
$$
    \int_{\bd\dd} \group{ \int_{\Lambda_{\zeta}^r} \abs{f}^q \dif\mu_{\alpha-1} }^{p/q} \dif l(\zeta) \simeq \int_{\bd\dd} \group{ \int_{\Gamma_{\zeta}} \abs{f}^q \dif\mu_{\alpha-1} }^{p/q} \dif l(\zeta)
$$
for any $r$. This phenomenon is called \emph{aperture invariance}. As such, the choice of $r = 1/\sqrt{2}$ is merely an arbitrary choice.

As mentioned before, we have normalized so that $T_{\alpha}^{p,p} = A_{\alpha}^p$. We may see this by applying Tonelli's theorem and noting that for $I_z = \set{ \zeta \in \bd\dd : z \in \Gamma_{\zeta} }$ we have
\[
    l(I_z) \simeq 1-\abs{z} \tag{$\star$}
\]
on $\dd$ where we write $f \simeq g$ on $E$ to mean
$$
    C^{-1} f(z) \leq g(z) \leq C f(z)
$$
for all $z \in E$ and $C>0$ is a constant independent of $z$. If the domain of $z$ is understood, we write $f \simeq g$. Note that we sometime take $z$ to denote a function and $f,g$ to denote semi-norms. We prove $(\star)$ now.

\Lemma{?1?}{
    Let $I_z = \set{\zeta \in \bd\dd : z \in \Lambda_{\zeta}}$ for some Stoltz family $\Lambda_{\zeta}$. Then $ l(I_z) \simeq 1-\abs{z}^2$.
}
\begin{proof} 
    By rotation invariance of $l(I_z)$, we can take $z = r$. Let $R$ be the aperture of $\Lambda_{\zeta}$. If $r \leq R$, then $l(I_z) = 1$. In this case, we choose $C = (1-R^2)^{-1}$ so that
    $$
        C^{-1} (1-\abs{z}^2) \leq 1 = l(I_z) \leq C(1-\abs{z}^2) \,.
    $$
    We now assume that $r>R$. We want to calculate the angle $\theta \in (0,\pi/2)$ such that $z \in \bd\Lambda_{e^{i\theta}}$. Equivalently by rotation invariance, we can calculate the angle $\theta\in (0,\pi/2)$ such that $w \defeq re^{i\theta} \in \bd\Lambda_{1}$. The upper boundary of the triangular portion of $\Lambda_{1}$ is given by the equation $y=M(1-x)$ for $M$ depending on $R$.
    Thus we can calculate
    $$
        w = \frac{ \group{ M^2 + \sqrt{ r^2(M^2+1)-M^2) } } + iM\group{ 1-\sqrt{r^2(M^2+1)-M^2}} }{ M^2 + 1}
    $$
    using quadratic formula. From this equation for $w$, the calculation follows directly.
\end{proof}

While we proved $(\star)$ directly, we may also compare with \citePage{6}{Perala}. In that paper (and others), the authors replace the Stoltz families with sets of the form
$$
    \bar{\Lambda}_{\zeta}^M = \set{ z \in \dd : \abs{1-z\bar{\zeta}} < M(1-\abs{z}) }
$$
or
$$
    \set{ z \in \dd : \abs{1-z\bar{\zeta}} < M(1-\abs{z}^2) }
$$
for $\zeta \in \bd\dd$ and $M>1/2$. In general, we define an \emph{approach family} to be a family of sets $E_{\zeta}$ for $\zeta \in \dd$. An \emph{approach region} is an element of a an approach family. Every approach family \emph{generates} a tent space by replacing $\Lambda_{\zeta}$ with $E_{\zeta}$ in the definition of the tent norm. We will say that two approach families are equivalent if they generate the same tent spaces. The \emph{standard approach family} is defined to be $\Gamma_{\zeta}$. The tent space generated by the standard approach region may be called the \emph{standard tent space} (or just the tent space as above). The tent spaces defined by the above approach families are equivalent. We see this by applying the following lemma together with the equivalence of the families $\Lambda_{\zeta}^r$ for $r \in (0,1)$ and a calculation.\Lemma{?2?}{
    If $E_{\zeta}^-$ and $E_{\zeta}^+$ are two equivalent approach families and $E_{\zeta}$ is another approach family such that
    $$
        E_{\zeta}^- \subseteq E_{\zeta} \subseteq E_{\zeta}^+
    $$
    for each $\zeta \in \bd\dd$, then the family $E_{\zeta}$ is equivalent to the families $E_{\zeta}^-$ and $E_{\zeta}^+$.
}
\begin{proof}
    The result follows from the inequalities
    $$
        \int_{E_{\zeta}^-} \abs{f}^q \dif \mu_{\alpha-1} \leq \int_{E_{\zeta}} \abs{f}^q \dif \mu_{\alpha-1} \leq \int_{E_{\zeta}^+} \abs{f}^q \dif \mu_{\alpha-1}
    $$
    and the fact that
    $$
        \int_{\bd \dd} \group{ \int_{E_{\zeta}^+} \abs{f}^q \dif \mu_{\alpha-1} }^{p/q} \dif l(\zeta) \leq C \int_{\bd \dd} \group{ \int_{E_{\zeta}^-} \abs{f}^q \dif \mu_{\alpha-1} }^{p/q} \dif l(\zeta)
    $$
    by the equivalence of $E_{\zeta}^-$ and $E_{\zeta}^+$.
\end{proof}
We now record this equivalence in the following lemma. The following claimed containments are a calculation. We omit the details. (Please see these references for more information: \cite[Theorem 1.3]{Pavlovic}, \citePage{236}{Pelaez}, \citePage{25}{Perala}.)
\Lemma{?3?}{
    For any $r$, there are $M^+$ and $M^-$ such that $\bar{\Lambda}_{\zeta}^{M^-} \subseteq \Lambda_{\zeta}^r \subseteq \bar{\Lambda}_{\zeta}^{M^+}$. As the families $\bar{\Lambda}_{\zeta}^{M^-}$ and $\bar{\Lambda}_{\zeta}^{M^+}$ are equivalent, $\bar{\Lambda}_{\zeta}^M$ generates the standard tent space for all $M$.
}

The following appears in \citePage{66}{Arsenovic}. (Although, the author first learned about it from \cite{Perala}.) A version of this lemma for the upper half plane appears in \cite{Luecking_Embedding} which predates the other references. 
\SpecialLemma{?C?}{ 
    Let $s \in (0,\infty)$, $\lambda > \Max{1,1/s}$, and
    $$
    K(z,\zeta) = \frac{1-\abs{z}^2}{\abs{1-\overbar{\zeta}z}} \,.
    $$
    For any measure $\mu$ we have 
    $$
      \oneover{C} \int_{\bd \dd} \mu(\Gamma_{\zeta})^s \dif l(\zeta) \leq 
      \int_{\bd \dd} \group{\int_{\dd}  K(z,\zeta)^{\lambda} \dif\mu(z) }^s \dif l(\zeta) \leq 
      C \int_{\bd \dd} \mu(\Gamma_{\zeta})^s \dif l(\zeta)
    $$
    where $C>0$ is independent of $\mu$. Moreover, the tent norm is independent of aperture. 
}
\noindent Strictly speaking, Arsenovic proved the above result for approach families $\bar{\Lambda}_{\zeta}^M$ in place of $\Lambda_{\zeta}^M$. From this and \thref{?3?}, the above result holds as written. 

The above lemma plays a crucial role in our understanding of tent spaces. We give a notation for tent spaces defined by the above kernels: We define
$$
    \myTilde{\scriptA}^q_{\alpha,\lambda}[f;\mu] (\zeta) = \group{ \int_{\dd} K(z,\zeta)^{\lambda} \abs{f(z)}^q \dif\mu_{\alpha-1}(z)}^{1/q}
$$
and
$$
    \norm{f}_{\myTilde{T}_{\alpha,\lambda}^{p,q}(\mu)} = \group{ \int_{\bd\dd}[ \myTilde{\scriptA}_{\alpha,\lambda}^q [f;\mu] ]^p \dif l }^{1/p} \,.
$$
Finally, we let $\myTilde{T}_{\alpha,\lambda}^{p,q}(\mu)$ be the space of measurable functions with finite norm. By \thref{?C?}, $\myTilde{T}_{\alpha,\lambda}^{p,q}(\mu) = T_{\alpha}^{p,q}(\mu)$ with equivalence of norms when $\lambda > \Max{1,q/p}$ --- in which case, we may omit the dependence of $\lambda$ when referring to $\myTilde{T}_{\alpha,\lambda}^{p,q}(\mu)$. As before, we drop the notation for the measure when $\mu = A$. We define $\tilde{t}_{\alpha,\lambda}^{p,q}$ in the analogous way to $t_{\alpha}^{p,q}$.

Given the definition of the tent spaces, one may not be surprised by the calculation of its dual when $p$ and $q$ are greater than one. See \cite{Perala} for the following result. 
 \SpecialTheorem{?D?}{
	If $p,q > 1$ then $T_{\alpha}^{p,q}$ is a Banach Space. We have that $(T_{\alpha}^{p,q})^* = T_{\alpha}^{p',q'}$ and $(AT_{\alpha}^{p,q})^* = AT_{\alpha}^{p',q'}$ under the pairing 
    $$
		\inner{f,g}_{\alpha} = \int_{\dd} f\overbar{g} \dif A_{\alpha} 
    $$
    where $p'$ and $q'$ are dual exponents to $p$ and $q$ respectively.
}
Please see \cite{Perala} for the dual of $AT_{\alpha}^{p,q}$ in the cases $p \leq 1$ and $q \leq 1$. See \cite{Luecking_Embedding} for similar results in the upper-half space. We will not need these results. We now observe that $T_{\alpha}^{p,q}(\mu)$ is a quasi-Banach space even when $p$ or $q$ is less than 1.

\Lemma{?4?}{
    Let $s = \Min{p,q,1}$. Then 
    $$
        \norm{f}_1 \defeq \norm{f}_{T_{\alpha}^{p,q}(\mu)}^s 
        \andEq 
        \norm{f}_2 \defeq \norm{f}_{\myTilde{T}_{\alpha}^{p,q}(\mu)}^s
    $$
    are $s$-norms which define metrics via $d_j(f,g) = \norm{f-g}_j$ for each measure $\mu$ and $j=1,2$.
}
\begin{proof}
    We sketch the proof that $\norm{\Dott}_1$ is an $s$-norm. The proof for $\norm{\Dott}_2$ is similar. Set $\hat{s} = \Min{q,1}$. We only need to verify that $f \mapsto \group{\scriptA_{\alpha}^q[f;\mu]}^{\hat{s}}$
    and $f \mapsto \norm{ f }_{L^{p/\hat{s}}}^{s/\hat{s}}$ satisfy the triangle inequality by directly checking all valid combinations of $s$ and $\hat{s}$.
\end{proof}

Finally, we collect some known results for function theory in the unit disc. We start by defining a Carleson measure for $A_{\alpha}^p$ to be any measure $\mu$ such that the inclusion
$$
    A_{\alpha}^p \to A_{\alpha}^p(\mu)
$$
is bounded. We have the characterization of such measures in \citePage{166}{Zhu} as follows:
\SpecialTheorem{?E?}{
    Let $\mu$ be a positive, finite measure. The following are equivalent
    \begin{enumerate}
        \item $\mu$ is a Carleson measure for $A_{\alpha}^p$.
        \item $\mu(D(z,r)) \leq C_r A(D(z,r))$ for all $r\in (0,1)$.
        \item $\mu(D(z,r)) \leq C_r A(D(z,r))$ for some $r\in (0,1)$.
    \end{enumerate}
}
When $\mu = \nu^{z\0}$ for some sequence $z\0$ we have the following result which is essentially found in \citePage{70}{DurenSchuster}.
\SpecialTheorem{?F?} {
    Let $z\0$ be a sequence in $\dd$. For each $\delta>0$ and $w \in \dd$, define the counting functions
    $$
        N(z\0\SP,w,\delta) = \# \group{ z\0 \cap D(w,\delta) } \andEq N(z\0\SP,\delta) = \sup_w D(z\0\SP,w,\delta)
    $$
    where $\#E$ is the number of elements in $E$. The following are equivalent.
  \begin{enumerate}
  \item $z\0$ is a finite union of uniformly discrete sequences.
  \item $N(z\0\SP,\delta)$ is finite for some $\delta \in (0,1)$.
  \item $N(z\0\SP,\delta)$ if finite for each $\delta \in (0,1)$.
  \item The measure $\nu^{z\0}$ is a Carleson measure for $A_{\alpha}^p$.
  \end{enumerate}
}
We will call any sequence satisfying \thref{?F?} a \emph{Carleson sequence}. We now collect some standard estimates and an equation in the following two lemmas. These will be used throughout our discussions.
 
\SpecialLemma{?G?}{
  Let $\delta \in (0,1)$. For $\rho(u,v) < \delta$ and $z \in \dd$ we have the following 3 standard estimates, 
  \begin{align*}
      (1) &\tab \oneover{C} \abs{1-z\overbar{u}} \leq \abs{1-z\overbar{v}} \leq C \abs{1-z\overbar{u}} \\
      (2) &\tab \oneover{C} (1-\abs{u}^2) \leq (1-\abs{v}^2) \leq C (1-\abs{u}^2) \\
      (3) &\tab \oneover{C} (1-\abs{z}^2)^{2} \leq A\group[\big]{ D(z,\delta) } \leq C (1-\abs{z}^2)^2
  \end{align*}
  for some $C$ depending only on $\delta$. From these it follows that
  $$
    \sigma_{\alpha}( D(z,\delta) ) \simeq (1-\abs{z}^2)^{1+\alpha}
  $$
  with equivalence constant depending only on $\delta$.
}

\SpecialLemma{?H?}{
    The following identity holds for $\zeta, z \in \dd$
    $$
        \abs{\phi_{\zeta}'(z)} = \frac{1-\abs{\phi_{\zeta}(z)}^2}{1-\abs{z}^2} \,.
    $$
}

When $p \leq 1$, we lose H\"older's inequality. However, we do gain an inequality of the form
$$
    \abs*{ \sum_k a_k }^p \leq \sum_k \abs{a_k}^p \,.
$$
Hence, we hope to find a way to turn integrals of $\abs{f}$ into sums when $f$ is holomorphic. The correct device for this purpose is a $\delta$-lattice. Given some $\delta>0$ we define a \emph{$\delta$-lattice} to be a set $z\0$ a uniformly discrete sequence with $\delta(z\0) \leq \delta$ such that $\set{D(z_k,2\delta)}_k$ covers $\dd$. We also use lattices to prove a result analogus to \thref{?F?} for the tent spaces. Luckily, the same argument as in lemma 4.8 of \citePage{70}{Zhu} guarantees the existence of such latices.

\SpecialLemma{?I?}{
    For each $\delta \in (0,1/2)$, there is a $\delta$-lattice.
}
\noindent We will see the utility of lattices throughout this paper.

\section{Analytic Tent Spaces in the Disc}

We begin by presenting some specialized results for analytic tent spaces in the disc. Many of the results here hold in some form for Bergman spaces. We also want to see how many of these results hold for sequence spaces of tent type. The overall goals of the section are as follows:
\begin{itemize}
    \item Define a family of quasi-isometries between tent spaces and sequence spaces of tent type.
    \item Characterize when the evaluation operator $E$ from before is bounded.
    \item To every tent space, associate a Bergman space which share some key properties. In particular, the spaces will have the same interpolating and sampling sequences.
\end{itemize}
In the process, we will develop key results for the later sections. 

We start by defining a class of functions which are similar enough to holomorphic functions so that we can prove boundedness results while still containing certain discontinuous functions. Let $M:(0,1/2) \arrow (0,\infty)$ be any function. We define $\scriptH_M$ to be the set of functions $f: \dd \arrow \cc$ with $\abs{f}$ being upper-semicontinuous such that for each $\delta \in (0,1/2)$ we have
$$
    \abs{f(z)} (1-\abs{z}^2)^{2} \leq M(\delta) \int_{D(z,\delta)} \abs{f(w)} \dif A(w) \,.
$$
Functions in $\scriptH_M$ are said to be quasi-subharmonic with defect $M$. We now exhibit some functions in $\scriptH_M$. The first example is $\abs{f}^p$ where $f$ is holomorphic.

\SpecialLemma{?J?}{
  If $f$ is holomorphic, then $\abs{f}^p$ is subharmonic. From this it follows that 
  $$
      \abs{f(z)}^p (1-\abs{z}^2)^{1+\alpha} \leq M(\delta) \int_{D(z,\delta)} \abs{f(w)}^p \dif\sigma_{\alpha}(w)
  $$
  for some $M(\delta) > 0$ depending only on $\delta$. Hence, we have that that $\holo(\dd) \subseteq \scriptH_M$.
}

The second class of functions are those of the form
$$
    F[a\0] = \sum_{k} a_k \chi_{D(z_k,\epsilon)}
$$
where $\epsilon>0$ is fixed and $a_k \in t_{\alpha}^{p,q}$. We can think of this mapping as an embedding of $t_{\alpha}^{p,q}$ into $T_{\alpha}^{p,q}$. The following 5 lemmas establish key features of the above embedding. In fact, the above embedding is quasi-isometric; although, it will take a few results to verify this. This embedding will allow us to translate many results for analytic tent spaces to sequence spaces of tent type (i.e. Proving results for the class $\scriptH_M \cap T_{\alpha}^{p,q}$ typically proves the result for both $AT_{\alpha}^{p,q}$ and $t_{\alpha}^{p,q}$.)

\Lemma{?5?}{
    Suppose $z\0$ is uniformly discrete and $\epsilon = \delta(z\0)$. Any function of the form 
    $$
        F[a\0] = \sum_{k=0}^\infty a_k \chi_{\overbar{D(z_k,\epsilon)}}
    $$
    is in $\scriptH_M$ where $M$ depends only on $\epsilon$. This also implies $\abs{F}^q \in \scriptH_M$.
}
\begin{proof}
    We will use $\abs{\dott}$ for Lebesgue measure in this proof. Take $w \in \dd$. If $w \notin \overbar{D(z_k,\epsilon)}$ for each $k$ then we are done. Assume that $w \in \overbar{D(z_k,\epsilon)}$ and let $z = z_k$. Then $\abs{F(w)} = \abs{a_k}$ and 
    $$
        \abs{a_k} \abs{ D(w,\delta) \cap D(z,\epsilon) } \leq \int_{D(w,\delta)} \abs{F(\zeta)} \dif A(\zeta) \,.
    $$
    Thus, we only need to bound $\abs{ D(w,\delta) \cap D(z,\epsilon) }$ from below by $M_{\epsilon}(\delta)(1-\abs{z}^2)$.
    
    Let $\dif\mu(\zeta) = (1-\abs{\zeta}^2)^{-2}\dif A(\zeta)$ be the M\"obius invariant measure. Using \thref{?G?} to estimate $1-\abs{\zeta}^2$ on $D(z,\epsilon)$, we obtain
    $$
        (1-\abs{z}^2)^2\mu\group[\big]{ D(z,\epsilon) \cap D(w,\delta) } \simeq \abs{ D(z,\epsilon) \cap D(w,\delta)}
    $$
    and 
    $$
        \mu\group[\big]{ D(0,\epsilon) \cap D(\phi_z(w),\delta) } \simeq \abs{ D(0,\epsilon) \cap D(\phi_z(w),\delta) }
    $$
    where the equivalence constants depend only on $\epsilon$. Combining these, we see that
    $$
        (1-\abs{z}^2)^{-2} \abs{ D(z,\epsilon) \cap D(w,\delta)} \geq 
        C_{\epsilon} \abs{D(0,\epsilon) \cap D(\phi_w(z),\delta)} \geq
        C_{\epsilon} \abs{D(0,\epsilon) \cap D(\epsilon,\delta)}
    $$
    The result is complete taking $M(\delta) = C_{\epsilon} \abs{ D(0,\epsilon) \cap D(\epsilon,\delta)}$.
\end{proof}

We now hope to show this is a quasi-isometry from $t_{\alpha}^{p,q}$ into $T_{\alpha}^{p,q}$; however, the computation initially seems to depend on the way the discs $D(z_k,\epsilon)$ intersect the Stoltz regions $\Gamma_{\zeta}$. We show this is not the case. This will require showing that we can replace the Stoltz regions with their pseudohyperbolic neighborhoods. Unfortunately, the pseudohyperbolic neighborhood of a Stoltz region is no longer a Stoltz region. Denoting the pseodohyperbolic $\delta$-neighborhood of $E$ by $N_{\delta}(E)$, we have the following result:
\Proposition{?6?}{
    For any Stoltz family $\Lambda_{\zeta}$ of aperture $r$ and any $\delta \in (0,1)$, there is a Stoltz family $\Lambda_{\zeta}^+$ with $N_{\delta}(\Lambda_{\zeta}) \subseteq \Lambda_{\zeta}^+$. Moreover there is a constant $K = K(r)$ such that whenever $\delta<K$ there is a Stoltz family $\Lambda_{\zeta}^-$  satisfying $N_{\delta}(\Lambda_{\zeta}) \subseteq \Lambda_{\zeta}^+$.
}

The above is a tedious calculation. We leave the proof to the reader. Combining \thref{?6?} with \thref{?2?} and aperture invariance, we get the following.
\Corollary{?7?}{
    The approach family $N_{\delta}(\Lambda_{\zeta}^r)$ generates the standard tent space for any $r \in\nobreak (0,1)$ and $\delta \in (0,\delta_r)$. More explicitly, there is a constant $C$, depending on $r$ and $\delta$, such that
    $$
        C^{-1} \norm{f}_{T_{\alpha}^{p,q}}^p \leq \int_{\bd\dd} \group{ \int_{N_{\delta}(\Lambda_{\zeta}^r)} \abs{f}^q \dif \sigma_{\alpha}}^{p/q} \dif l(\zeta) \leq C \norm{f}_{T_{\alpha}^{p,q}}^p
    $$
    for each $f \in T_{\alpha}^{p,q}$.
}

We are finally able to conclude that the map is a quasi-isometry.

\Proposition{?8?}{
    Let be $z\0$ any uniformly discrete sequence in $\dd$ and $\delta \leq \delta(z\0)/2$ sufficiently small. 
    Given $a\0 \in t_{\alpha}^{p,q}(z\0)$, define 
    $$
    F[a\0] = \sum_{n=1}^{\infty} a_n \chi_{\overbar{D(z_n,\delta)}} \,.
    $$ 
    Then $\norm{F[a\0]}_{T_{\alpha}^{p,q}} \simeq \norm{a\0}_{t_{\alpha}^{p,q}(z\0)}$.
}
\begin{proof}
    Let $F = F[a\0]$. We assume $\delta$ is small enough to apply \thref{?5?,?6?}. By \thref{?6?}, we can choose an aperture small enough that 
    $\Gamma_{\zeta}^- \defeq \Gamma_{\zeta}^M$ satisfies $N_{\delta}(\Gamma_{\zeta}^{-}) \subseteq \Gamma_{\zeta}$. By aperture invariance of the tent norm (\thref{?C?}) and Tonelli's theorem, we have
    $$
        \norm{F}_{T_{\alpha}^{p,q}}^p  \leq 
        C\int_{\bd \dd} \group{ \int_{\Gamma_{\zeta}^-} \sum_k \abs{a_k}^q \chi_{D(z_k,\delta)}  \dif\sigma_{\alpha} }^{p/q} \dif l(\zeta)  \leq
        C\int_{\bd \dd} \group{  \sum_{z_k\in \Gamma_{\zeta} } \abs{a_k}^q  \sigma_{\alpha}\group[\big]{ D(z_k,\delta) \cap \Gamma^-_{\zeta} } }^{p/q} \dif l(\zeta)
    $$
    where the second inequality follows from the fact that $\Gamma_{\zeta}^- \cap D(z_k,\delta) = \emptyset$ if $z_k \notin \Gamma_{\zeta}$\,. Using \thref{?G?} to estimate $\sigma_{\alpha}( D(z_k,\delta) )$, we obtain that
     $ \norm{F}_{T_{\alpha}^{p,q}} \leq C\norm{a\0}_{t_{\alpha}^{p,q}(z\0)}$.
    
    We again use \thref{?6?} to choose an aperture $N$ large enough that 
    $\Gamma_{\zeta}^+ \defeq \Gamma_{\zeta}^N$ satisfies $N_{\delta}(\Gamma_{\zeta}) \subseteq \Gamma_{\zeta}^+$. Again by aperture invariance (\thref{?C?}) and Tonelli's theorem, we get 
    $$
        \norm{F}_{T_{\alpha}^{p,q}}^p \geq 
        C\int_{\bd \dd} \group{ \int_{\Gamma_{\zeta}^+} \sum_k \abs{a_k}^q \chi_{D(z_k,\delta)}  \dif\sigma_{\alpha} }^{p/q} \dif l(\zeta) \geq
         C\int_{\bd \dd} \group{ \sum_{z_k \in \Gamma_{\zeta} } \abs{a_k}^q \sigma_{\alpha}\group[\big]{ {D(z_k,\delta) \cap \Gamma_{\zeta}^+} } }^{p/q} \dif l(\zeta)
    $$
    where the second inequality follows by removing $\set{z_k : z_k \notin \Gamma_{\zeta} }$ from the sum. As $z_k \in \Gamma_{\zeta}$ implies that $D(z_k, \delta) \subseteq \Gamma_{\zeta}^+$\,, we have that $D(z_k,\delta) \cap \Gamma_{\zeta}^+=D(z_k,\delta)$. Applying \thref{?G?} gives $\norm{F}_{T_{\alpha}^{p,q}}^p \geq C \norm{a\0}_{t_{\alpha}^{p,q}(z\0)}^p$.
\end{proof}

We now give an analogous result. This time, we do not require $z\0$ to be uniformly discrete.

\Corollary{?9?}{
    Let $z\0$ be Carleson in $\dd$ and $\delta$ be sufficiently small. Given $a\0 \in t_{\alpha}^{p,q}(z\0)$, define 
    $$
    F[a\0] = \sum_{n=1}^{\infty} a_n \chi_{\overbar{D(z_n,\delta)}}
    $$ 
    Then 
    $
        \norm{F[a\0]}_{T_{\alpha}^{p,q}} \leq C \norm{a\0}_{t_{\alpha}^{p,q}(z\0)} \,.
    $
    for $C$ independent of $a\0\SP$. Moreover, if $a\0$ is positive, then 
    $$
        \norm{a\0}_{t_{\alpha}^{p,q}(z\0)} \leq C'\norm{F[a\0]}_{T_{\alpha}^{p,q}}
    $$
    for $C'$ independent of $a\0\SP$.
}
\begin{proof}
    Let $F = F[a\0]$. As $z\0$ is carleson, we can write $z\0$ as a union of uniformly discrete sequences $z\0^1, \dots, z\0^N$. Let $\delta$ be small enough so we can apply \thref{?8?} to $z\0^j$ for all $j$. Write $a\0$ as the corresponding union $a\0^1, \dots, a\0^N$. Let $F_j = F[a\0^j]$ be the functions from \thref{?8?}. Then 
    $$
        \norm{F}_{T_{\alpha}^{p,q}} \leq \sum_j \norm{F_j}_{T_{\alpha}^{p,q}} \simeq \sum_j \norm{a\0^j}_{t_{\alpha}^{p,q}(z\0^n)} \simeq \norm{a\0}_{t_{\alpha}^{p,q}(z\0)}
    $$
    as we can think of $a\0^n$ as an element of $t_{\alpha}(z\0)$. Thus $\norm{F[a\0]}_{T_{\alpha}^{p,q}} \leq C \norm{a\0}_{t_{\alpha}^{p,q}(z\0)}$ as desired.
    
    Now assume $a\0$ is positive. Then $\norm{F}_{T_{\alpha}^{p,q}} \simeq \sum_j \norm{F_j}_{T_{\alpha}^{p,q}}$ as $F_j \leq F$ are positive functions. As $z\0^n$ is uniformly discrete, we get the reverse inequality.
\end{proof}

If the reader is familiar with the Bergman space theory of interpolation and sampling, the reader may recall maps of the form 
$f \mapsto (f\comp\phi_w) \phi_w'^{\gamma}$. For $\gamma = 2/p$ and $p \neq 2$, these maps are the only isometries of the Bergman space $A^p$. Furthermore, these maps play an integral part in characterizing sampling and interpolation sequences in the Bergman space; indeed, even the uniform densities are defined in terms of the maps $\phi_w$. This insight leads us to the following precise statement of M\"obius invariance in tent spaces.

\Theorem{?10?}{ 
    Let $\gamma(p,q) = 1/p + (1+\alpha)/q$. The operators
    $$
    	S_w^{p,q}f \defeq f\comp \phi_w \dott[\phi_w']^{\gamma(p,q)}
    $$
    are bounded from $T_{\alpha}^{p,q}$ to itself uniformly for $w \in \dd$. Moreover, the operators $S_w^{p,q}$ are isometries on the space $\myTilde{T}_{\alpha,\lambda}^{p,q}$ for $\lambda = 2q/p$. As $S_w^{p,q}$ is it's own inverse, the operators $S_w^{p,q}$ are uniformly bounded from below also.
}
\begin{proof} We first prove that $S_w = S_w^{p,q}$ isometry on the smoothed tent space $\myTilde{T}_{\alpha,\lambda}^{p,q}$ for $\lambda$ as given. 
    Take $s=p/q$ so that $\lambda = 2/s$ and define 
    $$
    	K(z,\zeta) = \frac{ 1-\abs{z}^2 }{ \abs{1-\overbar{\zeta}z } } \,.
    $$
    Use the change of variables $z \mapsto \phi_w(z)$ followed by $\zeta \mapsto \phi_w(\zeta)$ to get
    $$
    	\norm{S_wf}_{\myTilde{T}_{\alpha,\lambda}^{p,q}}^p = \int_{\bd \dd} \group{ \int_{\dd} K(\phi_w(z),\phi_w(\zeta))^{\lambda} \abs{f(z)}^q \abs{\phi_w'(z)}^{-1/s} \abs{\phi_w'(\zeta)}^{1/s} \dif \sigma_{\alpha}(z) }^{s}   \dif l(\zeta) \,.
    $$
    Note that we have used \thref{?H?} to simplify the Jacobians. We now calculate that 
    $$
    	K(\phi_w(z), \phi_w(\zeta) ) = K(z,\zeta) \abs{\phi_w'(z)}^{1/2} \abs{\phi_w'(\zeta)}^{-1/2} \,.
    $$
    For our special choice of $\lambda = 2/s$, we have
    $$
    	K(\phi_w(z),\phi_w(\zeta))^{\lambda} \abs{f(z)}^q \abs{\phi_w'(z)}^{-1/s} \abs{\phi_w'(\zeta)}^{1/s} = K(z,\zeta)^\lambda \abs{f(z)}^q \,.
    $$
    Thus $S_w$ is an isometry on $\myTilde{T}_{\alpha, \lambda}^{p,q}$, as desired.
    
    We must handle the remaining statement in three cases. The cases are $s<2$, $q>1$, and $q < 1$. 
    
    If $s < 2$, then $ \lambda > \Max{1,1/s}$. In this case, we apply \thref{?C?} to get,
    $$
        \norm{S_w^{p,q}f}_{T_{\alpha}^{p,q}} \leq 
        C \norm{S_w^{p,q}f}_{\myTilde{T}_{\alpha,\lambda}^{p,q}} =
        C\norm{f}_{\myTilde{T}_{\alpha,\lambda}^{p,q}} \leq 
        C' \norm{f}_{T_{\alpha}^{p,q}} 
    $$
    where $C$ does not depend on $w$. This gives the result in the case $s<2$. For the remainder of the proof, we assume that $s \geq 2$.
	
	Assume first that $q>1$. (As $s \geq 2$, this forces $p >2$.)	We now have that $(T_{\alpha}^{p,q})^* = T_{\alpha}^{p',q'}$ where $p'$ and $q'$ are dual exponents as in \thref{?D?}. 
	While tedious, we can check that $s \geq 1$ implies that $s' \defeq p'/q' \leq 1$. Now compute that $\abs{(S_w^{p,q})(f)} = \abs{ (S_w^{p',q'})^*(f) }$ via a change of variables and \thref{?H?}. Note also that $\phi_w'(\phi_w(z))\phi_w'(z) = 1$.
	Thus $\norm{S_w^{p,q}} = \norm{S_w^{p',q'}}$. As $s' < 2$, we can apply the above to get the result whenever $q > 1$. 
	
	If $q \leq 1$, then choose $Q>1$ and set $P = Q(p/q) \geq 2$. Define $R:T_{\alpha}^{p,q} \arrow T_{\alpha}^{P,Q}$ by $Rf = \abs{f}^{q/Q}$ so that $\norm{Rf}_{T_{\alpha}^{P,Q}}^P = \norm{f}_{T_{\alpha}^{p,q}}^p$ via direct calculation. Moreover, $\abs{S_w^{P,Q}Rf} = \abs{RS_w^{p,q}f}$ as 
	$$
	    (q/Q)\gamma(p,q) = (q/Q)(1/p + (1+\alpha)/q)) = 1/P + (1+\alpha)/Q = \gamma(P,Q) \,.
	$$
	But now $P,Q > 1$ so that we can apply the above to get 
	$$
	    \norm{S_w^{p,q}f}_{T_{\alpha}^{p,q}}^p =
	    \norm{RS_w^{p,q}f}_{T_{\alpha}^{P,Q}}^P = 
	    \norm{S_w^{P,Q}Rf}_{T_{\alpha}^{P,Q}}^P \leq 
	    C \norm{Rf}_{T_{\alpha}^{P,Q}}^P = 
	    C\norm{f}_{T_{\alpha}^{p,q}}^p
	$$
	implying that $S_w^{p,q}$ is uniformly bounded. This completes the cases.
\end{proof}

We now present the analogous result for sequence spaces. Note that that this is the most important application of \thref{?8?}.  

\Corollary{?11?}{
    Let $z\0$ be uniformly discrete and $w \in \dd$. Given $a\0 \in t_{\alpha}^{p,q}(z\0)$, define 
    $$
        z_n' = \phi_w(z_n) \andEq a_{n}' =\nobreak a_{n} \dott (\phi_w'(z_n))^{\gamma}
    $$
    for $\gamma = (1+\alpha)/q+1/p$. Then $a\0' \in t_{\alpha}^{p,q}(z\0')$ and $\norm{a\0}_{t_{\alpha}^{p,q}(z\0)} \simeq \norm{a\0'}_{t_{\alpha}^{p,q}(z\0')}$ with equivalence constant that does not depend on $z\0\SP$, $a\0\SP$, or $w$. Equivalently, the operators $S_w^{p,q}$ defined above are uniformly bounded from $T_{\alpha}^{p,q}(\nu^{z\0})$ to $T_{\alpha}^{p,q}(\nu^{\phi_w(z\0)})$.
}
\begin{proof} 
    Let $S = S_w^{p,q}$ and define
    $$
    	F = \sum_{n=1}^{\infty} a_n \chi_{D(z_n,\delta)}
    $$
    so that $\norm{a\0}_{t_{\alpha}^{p,q}(z\0)} \simeq \norm{F}_{T_{\alpha}^{p,q}}$ by \thref{?8?}. By \thref{?10?}, we have $\norm{F}_{T_{\alpha}^{p,q}} \simeq \norm{SF}_{T_{\alpha}^{p,q}}$
    with equivalence constant independent of $w$, $z\0\SP$, and $a\0\SP$. Applying the change of variables $z \mapsto \phi_w(z)$ and estimating $\abs{\phi_w'(z)}$ on $D(z_n,\delta)$ via \thref{?G?}, we obtain
    $$
        \abs{SF(z)} \simeq
        \sum_{n=1}^{\infty} \abs{a_n'} \chi_{D(z_n',\delta)}(z) \,.
    $$
    It follows that $\norm{SF}_{T_{\alpha}^{p,q}} \simeq \norm{a\0'}_{t_{\alpha}^{p,q}(z\0)}$ by \thref{?8?} again. This gives the desired norm equivalence with equivalence constant independent of $w$. The second statement follows directly from what we have done and the definition of $t_{\alpha}^{p,q}$.
\end{proof}

We have completed the first objective of this section. We now turn to the boundedness of the evaluation operator $E$ as defined above. For this, we want to characterize those measures $\mu$ such that the inclusion is bounded from $T_{\alpha}^{p,q}$ to $T_{\alpha}^{p,q}(\mu)$. Before we can do this, we need to do a calculation which is a variation on \thref{?1?}.

\Corollary{?12?}{
    Given $r > 0$ sufficiently small, define 
    $$
        I_z^{r} = \set{ \zeta \in \bd\dd : D(z,r) \subseteq \Gamma_{\zeta}} \,.
    $$
    Then $l(I_z^r) \simeq (1-\abs{z}^2)$.
}
\begin{proof}
    Take $I_z = \set{ \zeta \in \bd\dd : z \in \Gamma_{\zeta}}$. By definition, $I_z^r \subseteq I_z$ and thus
    $$
    \abs{ I_z^r } \leq \abs{I_z} \leq C(1-\abs{z}^2)
    $$
    by \thref{?1?}. Next, apply \thref{?6?} to get some $\Gamma_{\zeta}'$ and $R>0$ such that $z \in \Gamma_{\zeta}'$ implies $D(z,r) \subseteq \Gamma_{\zeta}$ for $r < R$. Finally,  apply \thref{?1?} to $J_z \defeq \set{\zeta \in \bd\dd : z \in \Gamma_{\zeta}' }$ to get $l(J_z) \simeq (1-\abs{z}^2)$ so that $J_z \subseteq I_z^r$ implies
    $$
    C'(1-\abs{z}^2) \leq l(J_z) \leq l(I_z^{r}) \,.
    $$
    Combining these inequalities gives the desired equivalence.
\end{proof}

With the calculation out of the way, we can proceed to the following theorem characterizing Carleson measures for the tent spaces:  

\Theorem{?13?}{
    Suppose $\mu$ is a positive, finite measure. The inclusion $AT_{\alpha}^{p,q} \to AT_{\alpha}^{p,q}(\mu)$ is bounded if and only if $\mu$ is a Carleson measure for some (or every) Bergman space. 
}
\begin{proof}
    We will denote the area of $E$ by $\abs{E}$ and define $\hat{\mu} = \mu_{\alpha-1}$. 
    We first assume that $\mu$ is a Carleson measure for some Bergman space. For $\delta>0$, fix a $(\delta/2)$-lattice $z\0\SP$.
    By \thref{?6?}, we can fix some Stoltz family $\Gamma_{\zeta}^+$ such that $N_{\delta}(\Gamma_{\zeta}) \subseteq \Gamma_{\zeta}^+$. 
    As $\set{D(z_k,\delta)}_{k=1}^{\infty}$ covers $\dd$, Tonelli's theorem implies
    $$
        \scriptA_{\alpha,\mu}^{q}f(\zeta) \leq
        \int_{\Gamma_{\zeta}}\sum_{k=1}^{\infty} \abs{f}^q \chi_{D(z_k,\delta)} \dif\hat{\mu} \leq 
        \sum_{z_k \in \Gamma_{\zeta^+} } \int_{D(z_k,\delta)} \abs{f}^q \dif\hat{\mu}
    $$
    as $z_k \notin \Gamma_{\zeta}^+$ implies $D(z_k,\delta) \cap \Gamma_{\zeta} = \emptyset$. Moreover, for any $z \in D(z_k,\delta)$, we have that 
    $$
        \abs{f(z)}^q (1-\abs{z}^2)^{\alpha+1} \leq 
        C \int_{D(z,\delta)} \abs{f}^q \dif\sigma_{\alpha} \leq 
        C\int_{D(z_k,2\delta)} \abs{f}^q \dif\sigma_{\alpha}
    $$
    as $D(z,\delta) \subseteq D(z_k,2\delta)$. Then applying the Carleson condition (\thref{?E?}) together with \thref{?G?} to estimate $\abs{f}^q$ on $D(z_k,\delta)$ we get
    $$
        \scriptA_{\alpha,\mu}^{q}f(\zeta) \leq 
        C \sum_{z_k \in \Gamma_{\zeta}^+} \frac{ \mu\group[\big]{D(z_k,2\delta)}}{(1-\abs{z_k}^2)^{2}} \int_{D(z_k,2\delta)} \abs{f}^q \dif\sigma \leq
        NC' \int_{\Gamma_{\zeta}^{++}} \abs{f}^q \dif\sigma_{\alpha}
    $$
    with $\Gamma_{\zeta}^{++}$ being a Stoltz angle chosen so that $N_{2\delta}(\Gamma_{\zeta}^+) \subseteq \Gamma_{\zeta}^{++}$ and $N$ is the maximum overlap of the discs $\set{D(z_k,2\delta)}$. (Such $N$ must exist by \thref{?F?}.) We now get that
    $$
        \scriptA_{\alpha,\mu}^{q}f(\zeta) \leq \int_{\Gamma_{\zeta}^+} \abs{f(z)}^q \dif\sigma_{\alpha}(z) \,.
    $$
    Applying $\norm{\dott}_{L^{p/q}(\bd\dd)}^q$ above and invoking aperture invariance (\thref{?C?}), yields the desired boundedness.
    
    Assume the inclusion $AT_{\alpha}^{p,q} \to AT_{\alpha}^{p,q}(\mu)$ is bounded. Define $f_w(z) = (1-\overbar{w}z)^{-2\gamma}$ for $\gamma = 1/p + (1+\alpha)/q$ and $w \in \dd$. Let $I_w^r$ be as in \thref{?12?} for $r$ small enough. On the disc $D(w,r)$, \thref{?G?} implies that
    $$
       (1-\abs{z}^2)^{-\gamma q} \simeq (1-\abs{w}^2)^{\gamma q} \abs{f_w(z)}^q = \abs{S_w(1)}^q
    $$
    where $S_w$ is as in \thref{?10?}. Hence, $\zeta \in I_w^r$ (or equivalently $D(w,r) \subseteq \Gamma_{\zeta}$) implies that
    $$
        ( 1 -\abs{w}^2 )^{-q/p} \frac{\mu(D(w,r))}{\abs{D(w,r)}} \leq
        C \int_{D(w,r)}(1-\abs{z}^2)^{-\gamma q} \dif \hat{\mu}(z) \leq
        C' \int_{\Gamma_{\zeta}} \abs{S_w(1)}^q \dif\hat{\mu} \,.
    $$
    Exponentiate the above by $p/q$, integrate over $I_w^r$, and apply \thref{?12?} to get 
    $$
        \groupFrac{\mu(D(w,r))}{\abs{D(w,r)}}^{p/q} \leq
        C \int_{I_w^r}  \group{ \int_{\Gamma_{\zeta}} \abs{S_w(1)}^q \dif\hat{\mu} }^{p/q} \dif l(\zeta) \leq 
        C \norm{S_w(1)}_{T_{\alpha}^{p,q}(\mu)}^p \,.
    $$
    Using the boundedness of the inclusion and \thref{?10?}, we get
    $$
        \frac{\mu(D(w,r))}{\abs{D(w,r)}} \leq C \norm{S_w(1)}_{T_{\alpha}^{p,q}}^{q} \leq C' \norm{1}_{T_{\alpha}^{p,q}}^{q} = C''
    $$
    concluding the result.
\end{proof}

The following corollary is a consequence of the characterization of Carleson sequences and the above:

\Corollary{?14?}{
    The operator $E: T_{\alpha}^{p,q} \arrow t_{\alpha}^{p,q}(z\0)$ is bounded if and only if $z\0$ is a Carleson sequence.
}

We have now been able to characterize boundedness of $E$. As computation of the tent norm is nontrivial, it is difficult to see if a  function belongs to $T_{\alpha}^{p,q}$. It turns out that for $\beta = (1+\alpha)p/q - 1$, we have that $A_{\beta-\epsilon}^p \subseteq AT_{\alpha}^{p,q}$ for each $\epsilon > 0$. Moreover, $AT_{\alpha}^{p,q} \subseteq A_{\beta+\epsilon}^p$ for each $\epsilon>0$. Thus in some sense, $AT_{\alpha}^{p,q}$ is close to $A_{\beta}^p$. Moreover, the same can be said about the sequence spaces $t_{\alpha}^{p,q}(z\0)$ and $\ell_{\beta}^{p}(z\0)$. As the critical density for interpolation and sampling on $A_{\beta}^p$ is $(\beta+1)/p = (\alpha+1)/q$, we have some indication that the main theorem holds. Indeed, for the proof of the interpolation theorem, we make critical use of these containments. We now verify these claims.

\Lemma{?15?}{
    Define $M_p(r,f) = \norm{f}_{L^p(\bd D(0,r), l)}$. If $f \in T_{\alpha}^{p,q}$ and $\abs{f}^q \in \scriptH_M$, then
    $$
      M_s(r,f) \leq C \norm{f}_{T_{\alpha}^{p,q}} \group{ {1-r^2} }^{-\bracket{(\alpha + 1)/q + 1/p - 1/s} }
    $$
    for each $s \in [p,\infty]$ and $C$ independent of $f$.
}
\begin{proof}
    As in theorem 5.9 of \citePage{84}{Duren}, we only need to check the cases $s=p$ and $s=\infty$. For completeness, we will recall this argument:
    For $q \in (p,\infty)$, we have
    $$
        M_q(r,f)^q = \int_{\bd\dd} \abs{f(r\zeta)}^{p} \abs{f(r\zeta)}^{q-p} \dif l(\zeta) \leq M_{\infty}(r,f)^{q-p} \dott M_p(r,f)^p \,.
    $$
    Exponentiation by $1/q$ and applying the assumption gives the result. We now verify the cases $s=p$ and $s=\infty$.
    
    We first take $s=p$. As $\set{r\zeta : r\in [0,1)} \subseteq \Gamma_{\zeta}^-$, \thref{?6?} lets us choose $\delta>0$ such that $D(r\zeta,\delta) \subseteq \Gamma_{\zeta}$ for each $\zeta$ and $r\geq 0$. Now estimate that 
    \[
        \abs{f(r\zeta)}^q(1-r^2)^{\alpha+1} \leq CM(\delta) \int_{D(r\zeta,\delta)} \abs{f}^q \dif\sigma_{\alpha} \leq CM(\delta) \int_{\Gamma_{\zeta}} \abs{f}^q \dif\sigma_{\alpha} \tag{$\star$}
    \]
    by \thref{?G?} and the definition of $\scriptH_M$. It follows that 
    $$
        M_p(r,f) \leq C\norm{f}_{T_{\alpha}^{p,q}} \group{ \oneover{1-r^2} }^{(\alpha + 1)/q }
    $$
    as desired
    
    We now take $s = \infty$. Let $\zeta \in I_z^\delta$ for $I_z^\delta$ as in \thref{?12?} and $\delta$ small enough. Then we have
    $$
        \abs{f(z)}^q(1-\abs{z}^2)^{\alpha+1} \leq CM(\delta) \int_{\Gamma_{\zeta}} \abs{f}^q \dif\sigma_{\alpha}
    $$
    by the same reasoning as $(\star)$. Exponentiating by $p/q$ and integrating over $I_{z}^{\delta}$ gives
    $$
        \abs{f(z)}^p(1-\abs{z}^2)^{(\alpha+1)p/q+1} \leq C \int_{I_z^{\delta}} \abs{f(z)}^p(1-\abs{z}^2)^{(\alpha+1)p/q} \dif l(\zeta) \leq  C'M(\delta) \norm{f}_{T_{\alpha}^{p,q}}^p
    $$
    by \thref{?12?}. Exponentiating by $1/p$, gives that
    $$
        \abs{f(z)} \leq C_{\delta}\norm{f}_{T_{\alpha}^{p,q}} (1-\abs{z}^2)^{-(\alpha+1)/q-1/p}
    $$
    As $(1-x^2)^{-(\alpha+1)/q-1/p}$ is increasing in $x$, the result follows.
\end{proof}

We note now the special case of $s = \infty$ in the below corollary. This case is crucial to understand how point evaluations interact with the norm of $f$. In particular, the map 
$$
    f \mapsto \sup_{\abs{z} \leq R < 1 } \abs{f(z)}
$$
is bounded in terms of $R$ so that any bounded family in $AT_{\alpha}^{p,q}$ is also a normal family. Thus we can apply Montel's Theorem to bounded families.

\Corollary{?16?}{
    If $f \in AT_{\alpha}^{p,q}$ and $\gamma = (1+\alpha)/q + 1/p$, then 
    \[
        \abs{f(z)}(1-\abs{z}^2)^{\gamma} \leq D \norm{f}_{T_{\alpha}^{p,q}}\,.  \tag{$\star$}
    \]
    As a result, point evaluations are uniformly bounded on compact sets and $AT_{\alpha}^{p,q}$ is complete.
}

The above result follows readily from \thref{?15?}, so we omit the details. We also have a similar result for the sequence space $t_{\alpha}^{p,q}(z\0)$ whenever $z\0$ is Carleson.

\Corollary{?17?}{
    If $a\0 \in t_{\alpha}^{p,q}(z\0)$ for a Carleson sequence $z\0$ in $\dd$ and $\gamma=(1+\alpha)/q + 1/p$, then 
    $$
        (1-\abs{z_k}^2)^{\gamma} \abs{a_k} \leq  C \norm{a_k}_{t_{\gamma}^{p,q}(z\0)}
    $$
    with $C$ independent of $a\0\SP$.
}
\begin{proof}
    Write $z\0$ to be the union of $z\0^1, \dots, z\0^N$. Assume $\delta < \delta(z\0^n)$ for all $n$. Set 
    $$
        F = F[\gap{\abs{a\0}}] = \sum_{k=1}^{\infty} a_k \chi_{D(z_k,\delta)}
    $$ 
    as in \thref{?5?,?9?}. Then $F \in \scriptH_M$ for some function $M$ so that \thref{?15?} applies with $s = \infty$. Let $\delta$ be small enough to apply \thref{?9?}. Then
    $$
        \abs{ a_k } \leq \abs{F(z_k)} \leq C \norm{F}_{T_{\alpha}^{p,q}}(1-\abs{z_k})^{-\gamma} \leq C' \norm{a\0}_{t_{\alpha}^{p,q}(z\0)}(1-\abs{z_k})^{-\gamma}
    $$
    so that multiplication yields the desired inequality.
\end{proof}

Before tackling the containment problem, we must present a simple lemma which will help to \enquote{replace} Holder's inequality throughout the paper in some cases when $p<1$, $q<1$ or $p/q<1$. 

\Lemma{?18?}{
    Suppose $z\0$ is a Carleson sequence and $z_n' \in D(z_n,\delta)$ for each $n$. 
    Then $z\0'$ is also a Carleson sequence.
}

This clearly follows from \thref{?F?}. We omit the details. With all the tools in place, we are ready to prove the following result which simultaneously addresses the case of sequence spaces and the case of area measure. 

\Proposition{?19?}{
    Set $\beta = (1+\alpha)p/q - 1$. Given a function $f$ such that $\abs{f}^q \in \scriptH_M$ we have
    $$
        \oneover{C_{\epsilon}}\norm{f}_{ L^p_{ \beta + \epsilon } } \leq 
        \norm{f}_{ T_{\alpha}^{p,q} } \leq C_{\epsilon}
        \norm{f}_{ L^p_{ \beta - \epsilon} }
    $$
    for $\epsilon > 0$. Note that $C_{\epsilon}$ may depend on $M$. Moreover, if $p \leq q$, then 
    $$
    \norm{f}_{T_{\alpha}^{p,q}} \leq \norm{f}_{ L^p_{\beta } }
    $$
    We note that $C_{\epsilon} \arrow \infty$ as $\epsilon \arrow 0$.
}
\begin{proof}
    Suppose $f \in T_{\alpha}^{p,q}$. We can apply \thref{?15?} to get
    $$
        M_p(f,r)^p \leq C \norm{f}_{T_{\alpha}^{p,q}}^p (1-\abs{r}^2)^{-(\beta+1)}
    $$ 
    which yields  
    $$
        \norm{f}_{L^p_{\beta + \epsilon}}^p \leq C \norm{f}_{T_{\alpha}^{p,q}}^p \int_0^1 (1-r^2)^{\epsilon-1} \dif r \leq C' \norm{f}_{T_{\alpha}^{p,q}}^p
    $$
    giving one of the desired ineqalities.
    
    Next assume that $s = p/q > 1$ and denote $s' = s/(s-1)$. For $\epsilon > 0$ and 
    $$
        \delta = (\beta - \epsilon)/s-\alpha = ((1+\alpha)s-1-\epsilon)/s - \alpha
    $$
    we have $s(\delta+\alpha)  = \beta - \epsilon$ and $-s'\delta > -1$. We apply Holder's inequality to $f \in L_{\beta - \epsilon}^p$ yielding 
    \begin{align*}
        \group[\bigg]{ \int_{\Gamma_{\zeta}} \abs{f}^q \dif\sigma_{\alpha} }^s &\leq  
        \group[\bigg]{ \int_{\Gamma_{\zeta}} (1-\abs{z}^2)^{-s'\delta} \dif\sigma }^{s/s'}
        \group[\bigg]{ \int_{\Gamma_{\zeta}} \abs{f}^p (1-\abs{z}^2)^{s(\alpha + \delta)} \dif\sigma }  \\&=
        C_{\epsilon} \group[\bigg]{ \int_{\Gamma_{\zeta}} \abs{f}^p (1-\abs{z}^2)^{\beta - \epsilon-1} \dif A }
    \end{align*}
    Integrating the above over $\bd\dd$ and using Tonelli's theorem, we conclude that 
    $$
        \norm{f}_{T_{\alpha}^{p,q}}^p
        \leq C_{\epsilon} \int_{\bd\dd} \abs{f(z)}^p (1-\abs{z}^2)^{\beta - \epsilon}  \int_{I_z} (1-\abs{z}^2)^{-1} \dif l(\zeta) \dif A(z) 
        = C_{\epsilon}' \norm{f}_{A_{\beta-\epsilon}^p}^p
    $$
    where the equality above follows from \thref{?1?}.
    Finally, assume that $s \leq 1$. Fix some $\delta$-lattice $z\0$ for $\delta$ small enough that \thref{?9?} applies and take  $f \in L_{\beta}^p$. For each $k$ choose $z_k' \in\nobreak \overbar{D(z_k,2\delta)}$ such that 
    $$
        a_k \defeq \abs{f(z_k')} = \Sup{ \abs{f(z)} : z \in D(z_k, 2\delta)}
    $$
    As $z\0$ is uniformly discrete, $z\0'$ is a Carleson sequence by \thref{?18?}. By \thref{?6?} we can choose a Stoltz family $\Gamma_{\zeta}^+$ such that $a \notin \Gamma_{\zeta}^+$ implies $D(a,2\delta) \cap \Gamma_{\zeta} = \emptyset$. Let 
    $$
        F = F[\gap{a\0}] = \sum_{k=0}^{\infty} a_k \chi_{\overbar{D(z_k,2\delta)}}(z)
    $$
    The definition of $a_k$ implies $\abs{f(z)} \leq F(z)$. Apply \thref{?9?} to get 
    $$
        \norm{f}_{T_{\alpha}^{p,q}}^p \leq 
        \norm{F}_{T_{\alpha}^{p,q}}^p \leq 
        C\int_{\bd\dd} \group{ \sum_{z_k \in \Gamma_{\zeta}} \abs{a_k}^q (1-\abs{z_k}^2)^{\alpha+1} }^s \dif l(\zeta) 
    $$
    Finally, we can pull $s\leq 1$ inside and apply \thref{?1?} with Fubini's theorem to obtain 
    $$
        \norm{f}_{T_{\alpha}^{p,q}}^q \leq
        C\sum_k a_k^p (1-\abs{z_k}^2)^{(\alpha+1)s} l(I_{z_k}) \leq
        C'\sum_k \abs{f(z_k')}^p (1-\abs{z_k}^2)^{(\alpha+1)s+1}  \leq
        C''' \norm{f}_{A_{\beta}^p}^p
    $$
    using \thref{?F?,?G?}. So that we complete the proof.
\end{proof}
 
Note that the hypothesis $\abs{f}^q \in \scriptH_M$ is only necessary if $q/p < 1$. We now state the following direct consequence of the above.
 
\Corollary{?20?}{
    Set $\beta = (1+\alpha)p/q - 1$. We have the following continuous inclusions
    $$
      A^p_{\beta - \epsilon } \subseteq AT_{\alpha}^{p,q} \subseteq A^p_{\beta + \epsilon}
    $$
    for $\epsilon > 0$. Moreover, if $p<q$, then $A^p_{\beta} \subseteq AT_{\alpha}^{p,q}$
}

With just a little work, we can also get the following result:

\Corollary{?21?}{
  Let $z\0$ be uniformly discrete and $\beta = (1+\alpha)p/q - 1$. Given $\epsilon > 0$ we have the following continuous inclusions
  $$
      \ell^p_{\beta - \epsilon}(z\0) \subseteq t_{\alpha}^{p,q}(z\0) \subseteq \ell^p_{\beta + \epsilon}(z\0) \,.
  $$
}
As the above result is not used in the subsequent, we omit the proof.

\section{The Key Operator}

In this section, we define an operator on $t_{\alpha}^{p,q}(z\0)$ by
$$
    a\0 \mapsto (1-\abs{z}^2)^{b} \sum_{k=1}^{\infty} \abs{a_k} \frac{ (1-\abs{w_k} )^{c+2} }{ \abs{1-\overbar{w_k}z}^{c+b + 2} }
$$
which we call $S_{b,c}[a\0]$ for some real numbers $b$ and $c$. (Note that this operator is well-studied in the Bergman space case.) We give conditions on the $b$ and $c$ under which this operator is bounded from $t_{\alpha}^{p,q}$ to $T_{\alpha}^{p,q}$. Namely, in \thref{?25?}, we will prove that the operator is bounded in the following cases
\begin{itemize}
    \item $-b < (1+\alpha)/q$ and $c$ is sufficiently large for arbitrary $p,q$.
    \item $-b < (1+\alpha)/q<c+1$ for $p,q > 1$. 
\end{itemize}
This will be enough to complete the sufficiency of the sampling condition and of the interpolation condition. The computations in this section are straightforward applications of Holder's inequality, standard estimates for well-known kernels, discreteization of analytic functions, and standard tent space ideas. In fact, many of the relevant ideas of the proof come from the Bergman space version of this operator. We also want to recognize proposition 6 of \citePage{14}{Perala}. The case where $q>1$ is essentially a generalization of their ideas.

As operators between the same space are easier at times, we work with a related operator 
$$
    f(z) \mapsto (1-\abs{z}^2)^b \int_{\dd} f(w) \frac{ (1-\abs{w}^2)^c }{ \abs{1-\bar{w}z}^{2+b+c} } \dif A(w)
$$
from $T_{\alpha}^{p,q}$ to itself and then use \thref{?8?} to translate to the desired operator on sequences. We must break the proof of boundedness into two cases: $q>1$ and $q \leq 1$. We now present a calculation from \cite{Perala} which we use shortly.
\SpecialLemma{?K?}{ 
    Suppose $r+t > s + 2 > \Max{r,t} \geq \Min{r,t} > 0$ and $s>-1$. For some $C>0$ independent of $\zeta$ and $z$ we have 
    $$
    \int_{\dd} \frac{ (1-\abs{w}^2)^s }{\abs{1-\bar{w}\zeta}^r\abs{1-\bar{w}z}^t} dA(w) \leq C \abs{1-\bar{\zeta}z}^{2+s-t-r} \,.
    $$
}
We now define an intermediate operator which we will use in the subsequent calculations. We show this operator is bounded for appropriate choices of parameters.

\Lemma{?22?}{
    Define the nonlinear operator $U_{b,c} : T_{\alpha}^{p,q} \arrow T_{\alpha}^{p,q}$ by 
    $$
        U_{b,c}f(z) = \group{ (1-\abs{z}^2)^{b}\int_{\dd} \abs{f(w)}^q \frac{(1-\abs{w}^2)^c}{\abs{1-\overbar{w}z}^{b+c+2} \dif A(w) }  }^{1/q}
    $$
    if $c+1 > \alpha + \Max{1,q/p}$, and $b > -( \alpha+1)$, then $U_{b,c}$ is bounded.
}
\begin{proof}    Let $U=U_{b,c}$ and choose $\lambda > \Max{1,q/p}$ large enough such that for 
    \begin{align*}
        s &=\lambda + b + \alpha  - 1 \\ 
        r &= \lambda \\
        t &= 2+b+c
    \end{align*} 
    we have $s+2 > t > 1$. As $c$, $b$, and $\lambda$ are sufficently large, the hypotheses for \thref{?K?} are satisfied. We now apply Tonelli's theorem followed by \thref{?K?} to the above get 
    \begin{align*}
    	\group[\big]{\myTilde{\scriptA}_{\alpha,\lambda}^q[Uf](\zeta)}^q & =
        C  \int_{\dd} (1-\abs{w}^2)^{c}\abs{f(w)}^q \int_{\dd} \frac{ (1-\abs{z}^2)^{s} }{\abs{1-\overbar{\zeta}z}^{r} \abs{1-\overbar{w}z}^{t} } \dif A(z) \dif A(w) \\& \leq
        C' \int_{\dd} (1-\abs{w}^2)^{c}\abs{f(w)}^q \oneover{ \abs{1-\overbar{\zeta}w}^{c -\alpha +1} } \dif A(w) \\& =
        C' [\myTilde{\scriptA}_{\alpha,\lambda'}^q f(\zeta)]^q
    \end{align*}
    for $\lambda' = c + 1 - \alpha$. Finally, $\lambda' > \Max{1,q/p}$ so that
    $$
        \norm{Uf}_{T_{\alpha}^{p,q}} \leq 
        C \norm{Uf}_{\myTilde{T}_{\alpha,\lambda}^{p,q}} \leq
        C''\norm{f}_{\myTilde{T}_{\alpha,\lambda'}^{p,q}} \leq
        C'''\norm{f}_{T_{\alpha}^{p,q}}
    $$
    using \thref{?C?} twice. Thus $U$ is bounded. 
\end{proof}

We also need a standard estimate for a particular integral. The following estimate is in \citePage{32}{DurenSchuster}:
\SpecialLemma{?L?}{
    Let $1<2+t<s$, then 
    $$
      \int_{\dd} \frac{ (1-\abs{w}^2)^t }{ \abs{1-\overbar{w}z}^s } \dif A(w) \leq C (1-\abs{z}^2)^{2+t-s}
    $$
    where $C>0$ is independent of $z$.
}

Using the above lemmas, we are able to establish the boundedness of the above defined operator in the following two lemmas. The first lemma handles the case where $q \in (1,\infty)$.

\Lemma{?23?}{ 
    Let $ 1 < q < \infty$. The operator 
    $$
      S_{b,c} f(z) \defeq (1-\abs{z}^2)^b \int_{\dd} f(w) \frac{ (1-\abs{w}^2)^c }{ \abs{1-\bar{w}z}^{2+b+c} } \dif A(w)
    $$
    is bounded from $T_{\alpha}^{p,q}$ to itself whenever $-qb < \alpha + 1 < q(c+1)$ and one of the following occur: 
    \begin{enumerate}
    \item $p > 1$
    \item $q(c+1)-\alpha > q/p$ and $b+c+2 > q/p$.
    \end{enumerate}
    Note that $S_{b,c}$ is defined on the closure of $L^1$ in $T_{\alpha}^{p,q}$ for $p<1$ or $q<1$.
}
\begin{proof} 
    We break the proof of (1) into two cases: $q \leq p$  and $p>1$. Set $S = S_{b,c}$ for fixed $b$ and $c$ and let $p'$ and $q'$ be the dual exponents to $p$ and $q$ for $p>1$ and $q>1$.
    Using the hypothesis $-qb < \alpha + 1 < q(c+1)$, we can fix some
    $$
    \epsilon \in \group{-\frac{c+1}{q'}, \frac{b}{q'}} \bigcap \group{ -\frac{b+\alpha+1}{q}, \frac{c-\alpha}{q} }
    $$
    so that $1 < 2 + b +\epsilon q'< 2+c+b$. Apply Holder's inequality followed by \thref{?L?} to get 
    \begin{align*}
    	\abs{Sf(z)}^q &\leq  
        (1-\abs{z}^2)^{b q} \group{ \int_{\dd} \frac{(1-\abs{w}^2)^{c+\epsilon q'}}{\abs{1-\overbar{w}z}^{2+b+c}} \dif A(w) }^{q/q'} 
        \group{ \int_{\dd}\abs{f(w)}^q \frac{(1-\abs{w}^2)^{c-\epsilon q}}{\abs{1-\overbar{w}z}^{2+b+c}} \dif A(w) } \\& 
        \leq (1-\abs{z}^2)^{bq}\group{ C (1-\abs{z}^2)^{\epsilon q'-b} }^{q/q'} 
        \group{ \int_{\dd}\abs{f(w)}^q \frac{(1-\abs{w}^2)^{c-\epsilon q}}{\abs{1-\overbar{w}z}^{2+b+c}} \dif A(w) }  \\&
        = C^{q-1} (1-\abs{z}^2)^{b+\epsilon q} \int_{\dd}\abs{f(w)}^q \frac{(1-\abs{w}^2)^{c-\epsilon q}}{\abs{1-\overbar{w}z}^{2+b+c}} \dif A(w) \,.
    \end{align*}
    Let $\hat{b} = b+\epsilon q$ and $\hat{c} = c - \epsilon q$. For $U = U_{\hat{b},\hat{c}}$, we have $\abs{Sf} \leq C \abs{Uf}$.
    
    By our choice of $\epsilon$, we have $\hat{b} > \alpha+1$ and $\hat{c}>\alpha$.
    Now if $q \leq p$, then we have $\hat{c}+1-\alpha > 1 \geq q/p$. Thus $U$ is bounded by \thref{?22?}. Finally we have
    $$
        \norm{Sf}_{T_{\alpha}^{p,q}} \leq C\norm{Uf}_{T_{\alpha}^{p,q}} \leq C' \norm{f}_{T_{\alpha}^{p,q}}
    $$
    completing the proof when $q \leq p$.
    
    Next if $q > p > 1$ then we use duality, \thref{?D?}, to get the result as $q/p > 1$ implies $q'/p' < 1$. Using Tonelli's theorem, we calculate as follows:
    $$
        \inner{ \abs{S_{b,c}f}, \abs{g}}_{\alpha} \leq 
        \int_{\dd} \int_{\dd} \frac{  (1-\abs{z}^2)^b\abs{g(z)}(1-\abs{w}^2)^c \abs{f(w)} }{ \abs{1-\overbar{w}z}^{2+b+c}}  \dif A(w) \dif A_{\alpha}(z) =
        \inner{ \abs{f}, S_{b',c'}\abs{g} }_{\alpha}
    $$
    where $b' = c-\alpha$ and $c' = b+\alpha$. We verify that $-q'b' < \alpha + 1 < q'(c'+1)$ so that $\inner{\abs{f},S_{b',c'}\abs{g}} < \infty$ by our above result. We can now use Fubini's theorem to calculate $S_{b,c}^* = S_{b',c'}$. Finally, $S_{b',c'}$ is bounded on $T_{\alpha}^{p',q'}$ by the above so that $S_{b,c}$ is bounded also on $T_{\alpha}^{p,q}$.
    
    Now if (1) fails, then we can use (2) to choose 
    $$
        \epsilon = \delta/q + \Max{-(c+1)/q', -(b+\alpha+1)/q }
    $$
    for $\delta > 0$. Moreover, $\delta$ can be taken sufficiently small so that
    $$
        \hat{c} + 1 > q/p > 1
    $$ 
    where $\hat{b} = b+\epsilon q$ and $\hat{c} = c - \epsilon q$. We apply \thref{?22?} again to conclude $U = U_{\hat{b},\hat{c}}$ is bounded. As above, $\abs{Sf} \leq C \abs{Uf}$ This concludes the boundedness result.
\end{proof}

Note that the above result holds when $q=1$ with a slight change in the proof. We will not use this fact, however. We further speculate that the boundedness of $S_{b,c}$ on $T_{\alpha}^{p,q}$ implies $-qb < \alpha + 1 < q(c+1)$ in the case $p,q \geq 1$. The above result may also be sharp in the general case. We also wonder if the following result is sharp:

\Lemma{?24?}{
    Let $ 0 < q \leq 1$. Given a function $f$ with $\abs{f}^q \in \scriptH_M$, the operator
    $S_{b,c}$ as above satisfies 
    $$
    \norm{S_{b,c}f}_{T_{\alpha}^{p,q}} \leq\nobreak C \norm{f}_{T_{\alpha}^{p,q}}
    $$
    where $C$ is independent of $f$ (but may depend on $M$) whenever $-qb < \alpha + 1$ and $q(c+2)>\alpha + 1 +\Max{1,q/p}$
}
\begin{proof}
    Let $S = S_{b,c}$ and choose a $\delta$-lattice $w\0\SP$. As $\abs{f}$ is upper-semicontinuous, we can choose $w_n' \in D(w_n,2\delta)$ such that 
    $$
    \Sup[\big]{\abs{f(z)} : z \in \overbar{D(w_n,2\delta)} } = \abs{f(w_n')}\,.
    $$
    We now apply \thref{?G?} on each element of the cover $\set{D(z_k,2\delta)}$ of $\dd$ to get
    $$
        \abs{Sf(z)}^q \leq 
        C^q (1-\abs{z}^2)^{bq} \group{ \sum_{n=1}^{\infty} \abs{f(w_n')} \frac{ (1-\abs{w_n'}^2)^{c+2}}{\abs{1-\overbar{w_n'}z}^{b+c+2}} }^q
    $$
    where $C$ depends only on $\delta$ and $M$. As $q\leq 1$, we can pull the $q$ inside the sum above to get
    \begin{align*}
    	\abs{Sf(z)}^q &
        \leq C(1-\abs{z}^2)^{bq} \sum_{n=1}^{\infty} \abs{f(w_n')}^q \frac{ (1-\abs{w_n}^2)^{q(c+2)} }{ \abs{1-\overbar{w_n}z}^{q(2+b+c)} } \\&
        \leq M(\delta)C'(1-\abs{z}^2)^{bq} \sum_{n=1}^{\infty} \int_{D(w_n',\delta)} \abs{f(w)}^q \frac{ (1-\abs{w}^2)^{q(c+2)-2} }{ \abs{1-\overbar{w}z}^{q(2+b+c)} } \dif A(w) \\&
        \leq M(\delta)C''(1-\abs{z}^2)^{bq} \int_{\dd} \abs{f(w)}^q \frac{ (1-\abs{w}^2)^{q(c+2)-2 } }{ \abs{1-\overbar{w}z}^{q(2+b+c)} } \dif A(w)
    \end{align*}
    using \thref{?F?} and the definition of $\scriptH_M$ followed by \thref{?G?}. 
    
    Let $\hat{b} = bq$ and $\hat{c} = q(c+2)-2$. For $U = U_{\hat{b},\hat{c}}$, the above gives $\abs{Sf} \leq C \abs{Uf}$. The assumptions guarantee that $\hat{b}$ and $\hat{c}$ are sufficiently large so that $U$ is bounded by \thref{?22?}. Hence $S$ is also bounded.
\end{proof}

We are now ready for the main result of this section. Although the following is a corollary of the above two lemmas, we wish to give it a name reflecting its importance throughout the remainder of our discussions.

\Proposition{?25?}{
    Let $w\0$ be a Carleson sequence. For real numbers $b$ and $c$ where $-b < (\alpha + 1)/q$ and $c$ is sufficiently large, the mapping $S_{b,c}:\nobreak t_{\alpha}^{p,q}(w\0) \arrow T_{\alpha}^{p,q}$ given by 
    $$
      S_{b,c}[a\0](z) = (1-\abs{z}^2)^{b} \sum_{k=1}^{\infty} a_k \frac{ (1-\abs{w_k} )^{ c + 2} }{ \abs{1-\overbar{w_k}z}^{c + b + 2} }
    $$
    is bounded, where the sum converges in $T_{\alpha}^{p,q}$. Moreover, if $p,q > 1$, then $c+1 \geq (1+\alpha)/q$ is large enough.
}
\begin{proof}
    Define $\hat{S}:T_{\alpha}^{p,q} \arrow T_{\alpha}^{p,q}$ by
    $$
        \hat{S} f(z) \defeq (1-\abs{z}^2)^{b} \int_{\dd} f(w) \frac{ (1-\abs{w}^2)^{c} }{ \abs{1-\bar{w}z}^{2+b+c} } \dif A(w)
    $$
    as in \thref{?23?,?24?} and set $S = S_{b,c}$. Let $S_N$ be the $N^{\text{th}}$ partial sum of $S$ and define 
    $$
        F_N\bracket{\gap{\abs{a\0}}} = \sum_{k=1}^{N} \abs{a_k} \chi_{\overbar{D(w_k,\delta)}}
    $$
    as in \thref{?9?} (with $\abs{a_k}$ replaced with zero for $k > N$). Then $F_N\bracket{\gap{\abs{a\0}}}$ is in the domain of $\hat{S}$ and is in $\scriptH_M$. By \thref{?G?} and Tonelli's Theorem, we get
    \begin{align*}
        \abs{S_N[a\0](z)} &\leq 
        (1-\abs{z}^2)^{b} \sum_{k=1}^{N} \abs{a_k} \frac{ (1-\abs{w_k} )^{ c + 2} }{ \abs{1-\overbar{w_k}z}^{c + b + 2} } \\&\leq
        C (1-\abs{z}^2)^{b} \sum_{k=1}^{N}  \int_{D(w_k,\delta)} \abs{a_k} \frac{ (1-\abs{w} )^{c} }{ \abs{1-\overbar{w}z}^{c + b + 2}} \dif A(w) \\&\leq
        C K \abs{\hat{S}F_N\bracket{\gap{\abs{a\0}}(z)}}
    \end{align*}
    as everything is positive. By \thref{?23?,?24?}, we can choose $c$ large enough so that $\hat{S}:T_{\alpha}^{p,q} \arrow T_{\alpha}^{p,q}$ is bounded so that \thref{?9?} implies
    $$
    \norm{ S_N[a\0] }_{T_{\alpha}^{p,q}} \leq CK \norm{ \hat{S} F_N\bracket{\gap{\abs{a\0}}} }_{T_{\alpha}^{p,q}} \leq C' \norm{ F_N\bracket{\gap{\abs{a\0}}} }_{T_{\alpha}^{p,q}} \leq  C'' \norm{a\0}_{t_{\alpha}^{p,q}(w\0)}
    $$
    which gives that $S_N$ converges in the topology of $T_{\alpha}^{p,q}$ to a bounded operator $S$. When $p,q>1$, we only need to take $c+1 \geq (1+\alpha)/q$ so that $\hat{S}$ is bounded in this case by \thref{?23?}.
\end{proof}

While we believe the above result is sharp for $p,q \geq 1$, it is possible that the boundedness of $S_{b,c}$ depends on the sequence $z\0$ when $p<1$ or $q<1$.

\section{Interpolation} 
We now begin our discussion of the interpolation theorem. Recall that $z\0$ is interpolating if the evaluation operator $E:AT_{\alpha}^{p,q} \arrow t_{\alpha}^{p,q}$ given by
$$
    Ef = f(z\0)
$$
is surjective. We now recall the main theorem of this section.

\MainTheorem{interpolationTheorem'}{
    A sequence $z\0 \subseteq \dd$ is interpolating for $AT_{\alpha}^{p,q}$ if and only if $z\0$ is uniformly discrete and $D^+(z\0) < (1+\alpha)/q$. 
}

\noindent We will start with the easy case, the sufficiency of the interpolation condition. We will then transition to the necessity of the interpolation condition. The necessity will require some extra tools developed for the case of Bergman spaces. In the end, we will consider a natural framework to generalize the concept of interpolation. This more general framework will justify our choice of $t_{\alpha}^{p,q}$ as the target space for the evaluation operator
\subsection{ Sufficiency of the Interpolation Condition }

Given the work of Seip, Duren, and Schuster (\citePage{220}{DurenSchuster}, \citePage{54}{Seip_Book}, and \citePage{720}{Seip_Survey}) in the Bergman spaces; we can prove the interpolation condition is sufficient by simply writing down a formula and verifying it works. The most difficult part of the problem is checking that the given formula defines a bounded operator. This was done in the previous section. The main result we use from the above authors is the following:

\SpecialTheorem{?M?}{
    Let $b>0$. If $z\0$ is uniformly discrete with $D^+(z\0) < b$, then there is a sequence $w\0$ and an analytic function $f$ such that $z\0 \cup w\0$ is uniformly discrete and 
    $$
        \abs{f(z)} \simeq \rho(z,z\0 \cup w\0)(1-\abs{z}^2)^{-b} \,.
    $$
}

Without further delay, we give the interpolation formula in the proof of the following theorem:

\Proposition{?26?}{
  Every uniformly discrete sequence $z\0$ with $D^+(z\0) < (1+\alpha)/q$ is interpolating for $T_{\alpha}^{p,q}$.
}

\begin{proof} 
    Using \thref{?M?}, we construct an analytic function $g$ and a sequence $z\0'$ containing $z\0$ such that $z\0'$ is uniformly discrete and 
    $$
    	\abs{g(w)} \simeq \rho(w,z\0') (1-\abs{w}^2)^{-b}
    $$
    for some $b \in (D^+(z\0), (1 + \alpha)/q)$. As subsequences of interpolating sequences are interpolating, we can assume that $z\0 = z\0'$. 
    We now follow the proof of theorem 1 from \citePage{ 162}{DurenSchuster}: Given $a\0 \in t_{\alpha}^{p,q}(z\0) $, define
    $$
    	f(z) = \sum_{k=1}^{\infty} a_k \groupFrac{g(z)}{g'(z_k)(z-z_k)} \groupFrac{1-\abs{z_k}^2}{1-\overbar{z_k}z}^{c+1} \,.
    $$
    We wish to show $f(z_k) = a_k$ and $f \in AT_{\alpha}^{p,q}$. 
    Denote the $n^{\text{th}}$ partial sums above by $s_n$ and the summands by $f_k$.
    For $z$ close enough to $z_k$, we have $\rho(z,z_k) = \rho(z,z\0)$. By the estimates from \cite{DurenSchuster}, we get
    $$
    	\abs{f_k(z)} \leq 
        C \abs{a_k} (1-\abs{z}^2)^{-b}   \frac{ (1-\abs{z_k}^2)^{b + c + 2} }{\abs{1-\overbar{z_k}z}^{c + 2} }
    $$
    where $C$ does not depend on $z$ or $k$. We choose $c$ such that 
    $$
        b + c  > (1 + \alpha)/q + 1/p 
    $$
    and apply \thref{?17?} to get
    $$
        (1-\abs{z_k}^2)^{b+c+2} \abs{a_k} \leq (1-\abs{z_k}^2)^{2}\norm{a\0}_{t_{\alpha}^{p,q}(z\0)} \,.
    $$
    Whenever $\abs{z} \leq R < 1$, we have that 
    $$
    	C(1-R)^{b+c+2} \abs{s_n(z)} \leq
        \sum_{k=1}^n (1-\abs{z_k})^{c + b + 2} \abs{a_k} \leq
        \norm{a\0}_{t_{\alpha}^{p,q}(z\0)} \sum_{k=1}^n (1-\abs{z_k}^2)^{2}\,.
    $$
    By \thref{?F?} applied to the constant function $1$, the series
    $$
        \sum_{k=1}^{\infty} (1-\abs{z_k}^2)^{2}
    $$
    converges; hence, the partial sums $s_n$ converge uniformly on compact subsets. It follows that $f$ is analytic and $f(z_k) = a_k$.  Define $\hat{b} = -b$, $\hat{c} = c + b$, and 
    $$
        S[a\0](z) = S_{\hat{b},\hat{c}}[a\0](z) = (1-\abs{z}^2)^{\hat{b}} \sum_k \frac{(1-\abs{z_k}^2)^{\hat{c}+2}}{\abs{1-\overbar{z_k}z}^{\hat{b}+\hat{c}+2}}\abs{a_k}
    $$
    as in \thref{?25?}. By the estimates above, we see that
    $$
    \norm{s_n}_{T_{\alpha}^{p,q}} \leq C \norm{ S\bracket{a\0} }_{T_{\alpha}^{p,q}} \leq C' \norm{a\0}_{t_{\alpha}^{p,q}(z\0) } < \infty
    $$
    so that $f \in AT_{\alpha}^{p,q}$ whenever $\hat{c}$ is sufficiently large as $-\hat{b} = b < (1+\alpha)/q$. We only need to choose $c$ large enough to make everything work. We have now seen that $f$ is the desired interpolating function.
\end{proof}

\subsection{ Necessity of the Interpolation Condition }
In contrast to the above, showing the the necessity of the interpolation condition is quite intricate. We proceed by showing that interpolation in a tent space implies interpolation in a related Bergman space. We will achieve this in three steps: First, we will prove that interpolating sequences are uniformly discrete. Second, we will introduce and discuss weak convergence of sequences in the disc. Finally, we will prove that a small perturbation of an interpolating sequence is still interpolating and conclude the result.
    
The first step is required to make the rest of the steps work. The second step is crucial. This step is the instrument we will use to prove that interpolation sequences in a tent space are interpolating for a Bergman space. The final step is necessary to get a strict inequality $D^+(z\0) < (1+\alpha)/q$. Without the final step, we can not rule out $D^+(z\0) = (1+\alpha)/q$.

\subsubsection{Uniform Discreteness}

To prove that interpolation sequences are uniformly discrete, we want a formula of the form 
$$
    \bigabs{\abs{f(z)}(1-\abs{z}^2)^{b} - \abs{f(w)}(1-\abs{w}^2)^{b}} \leq C \rho(z,w) \norm{f}_{T_{\alpha}^{p,q}}
$$
for some $b$. Via the interpolation condition, we will be able to construct the appropriate collection of functions to conclude that $z\0$ is uniformly discrete. Along the way to this formula, we collect some intermediate results which will find applications here and in the sampling section. The formula above will follow from the following two calculations which are inspired by \cite{DurenSchuster}:

\Lemma{?27?}{
    Let $r \in (0,1)$ and $b>0$. There is some $C>0$ depending on $r$ such that for any analytic function $f$ we have
    $$
    	(1-\abs{z}^2)^{b}\abs{f(w)-f(z) } \leq C \rho(z,w) \group{ \int_{D(z,r)} \abs{f}^q \dif\sigma_{qb - 1} }^{1/q}
    $$
    whenever $\rho(z,w) < r/4$.
}

    The above estimates follow in the same manner as in \citePage{199}{DurenSchuster}. We leave the proof to the interested reader. The above inequality, which has it's own application later in this section, has a left-hand side which is asymmetrical in $z$ and $w$. Although the above is be sufficient for our current purposes, later applications implore us to state the symmetric version of the above.

\Lemma{?28?}{
	Define $Sf(z) = (1-\abs{z}^2)^{b} f(z)$ for $b > 0$. For any $r \in (0,1)$, there is some $C>0$, depending on $r$, such that for any analytic function $f$ we have
    $$
    	\bigabs{ \abs{Sf(z)} - \abs{Sf(w)} } \leq C \rho(z,w) \group{ \int_{D(z,r)} \abs{f}^q \dif\sigma_{qb - 1} }^{1/q}
    $$
    whenever $\rho(z,w) < r/4$.
}
Again, these estimates follow similarly to \citePage{199}{DurenSchuster}. We omit the proof. The above will also be applied in the sampling section. We are now ready to derive the aforementioned estimate.

\Corollary{?29?}{
    Define $Sf(z) = (1-\abs{z}^2)^{\gamma} \abs{f(z)}$ for $\gamma = (1+\alpha)/q + 1/p$. There is some $C>0$ such that for any analytic function $f$ we have
    $$
    	\abs{Sf(z) - Sf(w)} \leq C \rho(z,w) \norm{f}_{T_{\alpha}^{p,q}}
    $$
    whenever $\rho(z,w) < r/4$ for $r$ sufficiently small.
}
\begin{proof}
    Take $r \in (0,1)$ and $I_z^{r} = \set{ \zeta \in \bd\dd : D(z,r) \subseteq \Gamma_{\zeta}}$ as in \thref{?12?}. Apply \thref{?G?} to estimate $1-\abs{\zeta}^2$ on $D(z,r)$ and \thref{?28?} with $b = \gamma$ to get
    $$
        \abs{Sf(z)-Sf(w)}^p \leq 
        C \rho(z,w)^p (1-\abs{z}^2) \group{ \int_{D(z,r)} \abs{f}^q \dif\sigma_{\alpha} }^{p/q} \leq 
        C \rho(z,w)^p (1-\abs{z}^2) \group{ \int_{\Gamma_{\zeta}} \abs{f}^q \dif\sigma_{\alpha} }^{p/q} \,.
    $$
    Now integrate both sides over $I_z^{r}$, divide by $1-\abs{z}^2$, and use \thref{?12?} to get
    $$
        \abs{Sf(z)-Sf(w)}^p \leq 
        C \rho(z,w)^p \int_{\bd\dd} \group{ \int_{\Gamma_{\zeta}} \abs{f}^q \dif\sigma_{\alpha} }^{p/q} \dif l(\zeta)
    $$
    by the definition of $I_z^r$. The result follows by taking $p^{\text{th}}$ roots.
\end{proof}
We now generalize the notion of interpolation constants from Bergman spaces to tent spaces.
\Lemma{?30?}{
	Given an interpolating sequence $z\0$ for $AT_{\alpha}^{p,q}$, there is a
    constant $M>0$ such that for any $a\0 \in t_{\alpha}^{p,q}(z\0)$ there exists 
    $f \in AT_{\alpha}^{p,q}$ with $f(z_n) = a_n$ and 
    $\norm{f}_{T_{\alpha}^{p,q}} \leq M \norm{a\0}_{t_{\alpha}^{p,q}(z\0)}$.
}

The above lemma follows readily from the closed graph theorem (alternatively the open mapping theorem) as in the case of Bergman spaces. The constant $M$ from the above is called the \emph{interpolation constant} and is denoted $M(z\0)$. We now come to the conclusion of this subsection. Without the following result, none of the subsequent results will hold. 
\Lemma{?31?}{
	If $z\0$ is interpolating for $T_{\alpha}^{p,q}$, then $z\0$ is uniformly discrete
}
    This proof is exactly the same as lemma 18 in \citePage{234}{DurenSchuster} where we replace their lemma 2 with our \thref{?29?}.

\subsubsection{Weak Convergence of Sequences }
Given the previous section, we primarily concern ourselves with uniformly discrete sequences. However, it seems natural to start this section in the framework of general sequences which are locally finite in the sense that $N(z\0\hspace{0.5pt},0,r) < \infty$ for all $r \in [0,1)$. We now state the following definition which can be thought of as a weak version of interpolation.
\Definition{zero_set}{
    We say a sequence $z\0$ is a \emph{zero set} for a space $X$ of analytic functions $f:\dd \arrow \cc$ if there is some nonzero $f \in X$ such that $f(z\0) = 0\0$ as a set with repetition and $f(z)$ is nonzero for $z \notin z\0\SP$. 
 }
We will now associate a set of sequences, $W(z\0)$, to $z\0$ such that $z\0$ interpolates for $A_{\beta}^p$ if and only if every element of $W(z\0)$ is a zero set for $A_{\beta}^p$. As we know $T_{\alpha}^{p,q} \subseteq A_{\beta+\epsilon}^p$ for $\beta = (1+\alpha)p/q -1$, we know zero sets for $T_{\alpha}^{p,q}$ are also zero sets for $A_{\beta+\epsilon}^p$ for $\epsilon>0$. We will then prove that interpolation sequences for $T_{\alpha+\epsilon}^{p,q}$ are interpolation sequences for $A_{\beta}^p$. We first define the following topology on sequences
\Definition{convergence_of_sequences}{
    Let $\set{z\0^n}_n$ be sequence of sequences. We say that $\set{z\0^n}_n$ is \emph{well-behaved} if the function $n \mapsto N(z\0^n,0,r)$ is eventually constant for almost every $r$. If $\set{z\0^n}_n$ is not well-behaved, we say it \emph{diverges weakly}. Assume for the remainder that $\set{z\0^n}_n$ is well-behaved. We say that $\set{z\0^n}_n$ \emph{converges weakly} to an infinite sequence $z\0$ if, for each $n$, there is a rearrangement $w\0^n$ of $z\0^n$ such that $w_k^n \arrow z_k$ as $n \arrow \infty$ for all $k$. We say that a sequence $\set{z\0^n}_n$ of sequences \emph{converges weakly} to a finite sequence $\set{z_k}_{k=1}^N$ if, for each $n$, there is a rearrangement $w\0^n$ of $z\0^n$ such that $w_k^n \arrow z_k$ as $n \arrow \infty$ for all $k=1, \dots, N$ and $\abs{w_k^n} \arrow 1$ for all other $k$.
}
It is clear that a sequence can \enquote{converge} to two different sequences; however, the above convergence is well defined when considering $z\0^n$ and $z\0$ as sets with repetition. 
\Definition{weak_limits}{
    Let $z\0$ be a sequence in $\dd$. We define $W(z\0)$ to be the closure of $\set{ \phi_w(z\0) : w \in \dd}$ in the topology of weak convergence of sequences.
}
Our definition is similar to the definition from \citePage{198}{DurenSchuster}; however, \cite{DurenSchuster} requires the rearrangements to be naturally ordered instead of well-behaved. Let us define this.
\Definition{natually_ordered}{
    We say that $z\0$ is \emph{naturally ordered} if the sequence $\abs{z\0}$ is increasing. 
}
In \citePage{212}{DurenSchuster}, a sequence $z\0^n$ is said to converge weakly to $z\0$ if and only if there are naturally ordered rearrangements $w\0^n$ and $w\0$ of $z\0^n$ and $z\0$ such that $w_k^n \arrow z_k$ for each $k$ as $n\arrow \infty$. This definition is not quite equivalent to ours. Consider the following sequences $z\0^n$ defined by 
$$
    z_1^n = 1/2 + (-1)^{n}/2^{3n} \andEq z_2^n = -1/2 + (-1)^{n}/2^{3n}
$$
with $z_k^n = i(1-2^{-4k})$ for $k \notin {1,2}$. For $z_1 = 1/2$, $z_2 = -1/2$ and $z_k = i(1-2^{-4k})$, we have $z_k^n \arrow z_k$ as $n \arrow \infty$. Thus $z\0^n \arrow z\0$ weakly in our definition as $\set{z\0^n}_n$ is well-behaved. 

However, $\abs{z_1^n} < \abs{z_2^n}$ when $n$ is odd and $\abs{z_1^n} > \abs{z_2^n}$ when $n$ is even. Thus the only naturally ordered rearrangement of $z\0^n$ is $w_k^n = z_k^n$ if $n$ is odd and $w_k^n = z_{\phi(k)}^n$ if $n$ is even where $\phi(1) = 2$, $\phi(2)=1$, and $\phi(k) = k$ for $k \notin \set{1,2}$. The sequences $w_1^n$ and $w_2^n$ do not converge. Thus $z\0^n$ does not converge in the definition from \cite{DurenSchuster}. However, any subsequence of $z\0^n$ has a further subsequence which converges to $z\0$ in the definition of \cite{DurenSchuster}. Thus the convergence from \cite{DurenSchuster} can not define a topology. 

As the mode of convergence from \cite{DurenSchuster} is not topological, we use this alternative definition. Luckily, the necessary results from \cite{DurenSchuster} still hold upon replacing \enquote{naturally ordered} with \enquote{well-behaved}. We will recall some results from \cite{DurenSchuster} without proof unless the proofs need some significant modifications.

We now refocus ourselves on uniformly discrete sequences. We first ask when is the limit of uniformly discrete sequences still uniformly discrete. This leads to the following definition.
\Definition{equidiscrete}{
    A sequence $z\0^n$ of sequences is \emph{equidiscrete} if there is some $\delta > 0$ such that $\delta(z\0^n) \geq \delta$ for all $n$.
}
\noindent We hope that the limit of an equidiscrete sequence of sequences is uniformly discrete. 

\SpecialLemma{?N?}{
    Every equidiscrete sequence $z\0^n$ has a subsequence $z\0^{n_k}$ that converges weakly to a uniformly discrete sequence.
}

See \citePage{212}{DurenSchuster} for the proof. Also, we recall lemma 12 from \citePage{213}{DurenSchuster}.

\SpecialLemma{?O?}{
    Suppose $z\0^n$ satisfies $z_k^n \arrow z_k$ for some sequence $z\0$. If $f^n$ is a sequence of analytic functions converging uniformly on compact subsets to $f$, then $f_n(z_k^n) \arrow f(z_k)$.
}

Armed with the above, we present the crucial generalization of lemma 21 in \citePage{239}{DurenSchuster}. As the proof is essentially the same (noting that things work after replacing \enquote{naturally ordered} with \enquote{well-behaved}), we leave it out.
\Proposition{?32?}{
	Let $p \in (0,\infty)$ and $\alpha \in (0, \infty)$. A uniformly discrete sequence $z\0$ is interpolating for $A_{\alpha}^{p}$ if and only if every element of $W(z\0)$ is a zero set for $A_{\alpha}^p$.
}

In light of this proposition, we would like to see that if $z\0$ interpolates for a tent space then every element of $W(z\0)$ is a zero set for the tent space. Knowing this, we will prove that every element of $W(z\0)$ is a zero set for certain Bergman spaces --- the ones from \thref{?20?}. Indeed, we will soon see that each element of $W(z\0)$ is an interpolation sequence for the tent space provided that $z\0$ is interpolating. As in the Bergman space, we need all elements of $\set{ \phi_w(z\0) : w \in \dd}$ to be interpolating. This follows from our version of M\"obius invariance. Even more is true:
\Lemma{?33?}{
	If $z\0$ is interpolating for $AT_{\alpha}^{p,q}$, then \mbox{$\sup_{w \in \dd} \set{M(\phi_w(z\0))} < \infty$}. In particular, $\phi_w(z\0)$ is interpolating for $AT_{\alpha}^{p,q}$ for each $w\in\dd$.
}
\begin{proof}
    For $w \in \dd$, set $z\0' = \phi_w(z\0)$, and define $\delta = \delta(z\0)/2$. Fix $a\0' \in t_{\alpha}^{p,q}(z\0')$ and define
    $$
        a_{n} = a_{n}' \dott (\phi_w'(z_n'))^{\gamma} = a_{n}' \dott (\phi_w'(z_n))^{-\gamma}
    $$
    where $\gamma = (1+\alpha)/q+1/p$. By \thref{?11?}, we conclude that $a\0 \in t_{\alpha}^{p,q}(z\0)$ and $\norm{a\0}_{t_{\alpha}^{p,q}(z\0)} \leq C \norm{a\0'}_{t_{\alpha}^{p,q}(z\0')}$ for $C$ depending only on $p$, $q$, $\delta$, and $\alpha$.  As $z\0$ is interpolating, there is $f \in T_{\alpha}^{p,q}$ such that 
    $$
    \norm{f}_{T_{\alpha}^{p,q}} \leq M(z\0) \norm{a\0}_{t_{\alpha}^{p,q}(z\0)} \leq C M(z\0) \norm{a\0'}_{t_{\alpha}^{p,q}(z\0)}
    $$
    and $f(z_n) = a_n$.
    Setting $F(z) = (f \comp \phi_w(z)) \phi_w'(z)^{\gamma}$ and applying \thref{?10?}, we get 
    $$
        F(z_n') = a_n \dott (\phi_w'(z_n))^{\gamma} = a_n'
    $$
    and 
    $$
    	\norm{F}_{T_{\alpha}^{p,q}} \leq C \norm{f}_{T_{\alpha}^{p,q}} \leq C M(z\0) \norm{a\0'}_{t_{\alpha}^{p,q}(z\0)}
    $$
    so that $M( z\0' ) \leq C M(z\0)$ as $a\0'$ was arbitrary.
\end{proof}

With the above result, we are able to show each element of $W(z\0)$ interpolates for the tent space. Note that it is crucial that the interpolation constants $\set{M(\phi_w(z\0)) : w \in \dd}$ have an upper bound. We should not expect weak limits $z\0^n \arrow z\0$ where $M(z\0^n) \arrow \infty$ to have $z\0$ interpolating. Our result about $W(z\0)$ is a consequence of the following:

\Lemma{?34?}{
	Let $\set{z\0^n}_{n=0}^{\infty}$ be an equidiscrete sequence of sequences. If $\set{z\0^n} \arrow z\0$ weakly and $z\0^n$ is interpolating for $AT_{\alpha}^{p,q}$ for each $n$ with $\sup_n\set{M(z\0^n)} \leq M < \infty$, then $z\0$ is interpolating for $AT_{\alpha}^{p,q}$. Moreover, every element of $W(z\0)$ is interpolating if $z\0$ is interpolating.
}
\begin{proof}
	If $z\0$ is finite, we can interpolate with a polynomial; thus, we assume $z\0$ is infinite. 
	By rearranging, we assume that $z_k^n \arrow z_k$ as $n \arrow \infty$ for each $k$. Let $4\delta$ be the equidiscrete constant and $\gamma = (1+\alpha)/q+1/p$. 
	
	For each $N$, we can find $n_N$ such that 
	$$
	    \rho( z_k^n, z_k ) < \delta
	$$
	for $k \leq N$ and $n \geq n_N$. Replace $\set{z\0^n}_n$ with the subsequence $\set{z\0^{n_N}}_{N}$ which still converges to $z\0\SP$. Now for $n \leq k$, set 
	$$
	    a_k^n = a_k \groupFrac{ 1-\abs{z_k}^2 }{ 1-\abs{z_k^n}^2 }^{(1+\alpha)/q }
	$$
	and set $a_k^n = 0$ if $n > k$. Using \thref{?6?}, we choose $\Gamma_{\zeta}^+$ such that $z \in \Gamma_{\zeta}$ implies $D(z,\delta) \subseteq \Gamma_{\zeta}^+$. Thus $z_k^n \in \Gamma_{\zeta}$ implies $z_k \in \Gamma_{\zeta}^+$ for $k \leq n$. With \thref{?C?}, we now compute
	$$
	    \norm{a\0^n}_{t_{\alpha}^{p,q}(z\0^n)}^p \leq
	    \int_{\bd\dd} \group{ \sum_{k=0}^n \abs{a_k}^q (1-\abs{z_k}^2)^{1+\alpha} \chi_{\Gamma_{\zeta}^+}(z_k)}^{p/q} \dif \zeta  \leq
	    C \norm{a\0}_{t_{\alpha}^{p,q}(z\0)}
	$$
	so that $a\0^n \in t_{\alpha}^{p,q}(z\0^n)$ with uniform upper bound on the norms.
	
	By the interpolation hypothesis, we can find $f^n \in T_{\alpha}^{p,q}$ such that $f^n(z_k^n) = a_k^n$ and
    $$
    \norm{f^n}_{T_{\alpha}^{p,q}} \leq C M(z\0^n) \norm{a\0^n}_{t_{\alpha}^{p,q}(z\0)} \leq C' M \norm{a\0}_{t_{\alpha}^{p,q}(z\0)}
    $$
    since $M(z\0^n) \leq M$. By \thref{?16?}, we see that $\set{f^n}$ is a normal family. By Montel's theorem, we extract a subsequence of $\set{f^n} \subseteq AT_{\alpha}^{p,q}$ which converges to the function $F \in T_{\alpha}^{p,q}$ uniformly on compact sets. As $a_k^n \arrow a_k$ for each $k$ as $n \arrow \infty$, we conclude $F(z_k) = a_k$ by \thref{?O?} implying that $z\0$ is interpolating.
    
    Finally, suppose that $w\0 \in W(z\0)$. By definition, there is $\zeta\0$ such that for $w\0^k = \phi_{\zeta_k}(z\0)$ we have $w\0^k \arrow w\0$ weakly. It is clear that $\delta(w\0^k) = \delta(w\0)$ so that $\set{w\0^k}_k$ is equidiscrete. Applying \thref{?33?}, we get an upper bound on $\set{M(w\0^k)}_k$ so that $w\0$ is interpolating as desired.
\end{proof}

\subsubsection{Perturbation and the Conclusion}
We saw that any interpolation sequence $z\0$ for $T_{\alpha}^{p,q}$ interpolates for $A_{\beta + \epsilon}^p$ for $\beta = (1+\alpha)p/q - 1$. However, this is only good enough to prove $D^+(z\0) \leq (1+\alpha)/q$. We must apply this result to a suitable modification of $z\0\SP$. For each $a \in \dd$ we construct auxiliary sequences $z\0^a$ that is still interpolating which satisfies 
$$
    D(z\0\SP,a,r) \leq (1-\eta)D(z\0^a,a,r) + \epsilon(r)
$$
where $\epsilon(r) \arrow 0$ as $r \arrow 1$ and $\eta>0$ is some sufficiently small constant.
 
We soon define a function $g_w^{\eta}$ which we will use to form these auxiliary sequences. We will need to know how $g_w^{\eta}$ changes pseudohyperbolic distances and densities. We record the details of this function in the following technical lemma which is a modification and condensation of the content in page 203 of \cite{DurenSchuster}. 
 
\SpecialLemma{?P?}{
	Given $w \in \dd$ and $\eta \in (-1,1)$, define $g_{w}^{\eta} = \phi_{w}\comp g^{\eta} \comp\phi_{w}$ for 
    $$
    	g^{\eta}(z) = \groupFrac{\abs{z} + \eta }{1+\abs{z}\eta} \frac{z}{\abs{z}} \,.
    $$
    If $z\0$ is uniformly discrete, then the following hold:
    \begin{enumerate}
		\item $\rho( g_{w}^{\eta}(z), z) = \eta = \rho( g_{w}^{-\eta}(z),z )$ for each $w,z \in \dd$.
        \item For $\eta < \delta(z\0)/4$, we have $\delta( g_w^{\eta} (z\0) ) \geq \delta(z\0)/2$ and $\delta( g_w^{-\eta} (z\0) ) \geq \delta(z\0)/2$ for each $w \in \dd$.
        \item There is $C=C(\delta(z\0))>0$, such that for $\eta \in(0,1/4]$, $r \in (1/2,1)$ and any $w \in \dd$ we have 
        $$
			 D(z\0\SP,w,r) \geq (1-\eta)^{-1} \group{ D(g_w^{\eta}(z\0), w, r) - C\dott\bracket[\big]{\log(1-r)}^{-1} }
        $$
        and 
        $$
             D(z\0\SP,w,r) \leq (1-\eta) \group{ D(g_w^{-\eta}(z\0), w, r)  + C'\dott\bracket[\big]{\log(1-r)}^{-1} }
        $$
		where $(1+\eta) C' = (1-\eta)^{-1} C$.
	\end{enumerate}
}
We now discuss how the above can be obtained from \cite{DurenSchuster}: Part (1) follows from M\"obius invariance of the pseudohyperbolic metric together with \cite{DurenSchuster}. Part (2) follows directly from part (1). The first inequality of (3) is lemma 4 of \citePage{203}{DurenSchuster} combined with the fact that $D(z\0\SP,w,r) = D(\phi_w(z\0)\,, 0, r)$. The second inequality of (3) is obtained from the first by noting that $(g^\eta)^{-1} = g^{-\eta}$, 
$$
    g_{w}^{\eta} \comp g_{w}^{-\eta} = \phi_w \comp g^{\eta} \comp g^{-\eta} \comp \phi_w = \phi_w \comp \phi_w = \id
$$
and similarly that $g_w^{-\eta} \comp g_w^{\eta} = \id$. We then apply the first inequality with $z\0$ replaced with $g_w^{-\eta}(z\0)$.

One can think of the function $g_a^{\eta}$ as a uniform "push-towards-the-boundary" of points with the point $a$ thought of as the origin. Loosely speaking, the content of the above lemma is the following: Pushing a sequence towards the boundary decreases it's density, and pulling a sequence towards the center of the disc increases it's density.

We now show that interpolation is closed under small perturbations.

\Lemma{?35?}{
	Given any interpolation sequence $z\0$ for $AT_{\alpha}^{p,q}$ there is some $r > 0$ such that whenever $\rho(z_n,z_n') < r$ we have that $z\0'$ is also interpolating.
}
\begin{proof} 
    Roughly, we follow lemma 4 of \citePage{49}{Seip_Book}.
    Given $a\0 \in t_{\alpha}^{p,q}(z\0')$, let $a_k^0 = a_k$ so that $$
        \norm{a_k^0}_{t_{\alpha}^{p,q}(z\0)} \leq C \norm{a_k}_{t_{\alpha}^{p,q}(z\0')}
    $$
    by \thref{?G?,?6?}. We also applied aperture invariance, \thref{?C?}. 
    We can now find some $f^0 \in T_{\alpha}^{p,q}$ such that $f^0(z_k) = a_k^0$ and 
    $$
        \norm{f^0}_{T_{\alpha}^{p,q}} \leq M \norm{a\0^0}_{t_{\alpha}^{p,q}(z\0)}
    $$
    with $M = M(z\0)$\,. For any integer $n \geq 0$, assume that $f^n \in T_{\alpha}^{p,q}$ has been chosen such that $f^n(z_k) = a_k^n$ and 
    $$
        \norm{f^n}_{T_{\alpha}^{p,q}} \leq CM(\delta CM)^n \norm{a\0}_{t(z\0')} \,.
    $$ 
    Take $a_k^{n+1} = f^n(z_k) - f^n(z_k') $.
    By \thref{?33?}, we have $z\0$ is uniformly discrete. For $r < \delta(z\0)/4$\,, \thref{?27?} implies
    \begin{align*}
    	\norm{a\0^{n+1}}_{t_{\alpha}^{p,q}(z\0)}^p &=
        \int_{\bd \dd} \group[\Bigg]{\myspace \sum_{z_k \in \Gamma_{\zeta}} \bigabs{ (1-\abs{z_k})^{b} \bracket[\big]{ f^n(z_k) - f^n(z_k')} }^q }^{p/q} \dif l(\zeta) \\&\leq
        \int_{\bd \dd} \group[\Bigg]{ \myspace \sum_{z_k \in \Gamma_{\zeta}} [\delta C(r,q,\alpha) ]^q \group{ \int_{D(z_k,r)} \abs{f^n}^q \dif\sigma_{\alpha} } }^{p/q} \dif l(\zeta) \\& \leq
        [\delta C(r,q,\alpha)]^p \norm{f^n}_{T_{\alpha}^{p,q}}^p
    \end{align*}
    where the last inequality follows from \thref{?C?,?6?}. Thus, $a\0^{n+1} \in t_{\alpha}^{p,q}(z\0)$ so we can find $f^{n+1}\in T_{\alpha}^{p,q}$ such that $f_k^{n+1}(z) = a_k^{n+1}$ for all $k$ and
    $$
    	\norm{f^{n+1}}_{T_{\alpha}^{p,q}} \leq M \norm{a\0^{n+1}}_{t_{\alpha}^{p,q}(z\0)} \leq \delta C M \norm{f^n}_{T_{\alpha}^{p,q}} \leq CM(\delta CM)^{n+1} \norm{a\0}_{t_{\alpha}^{p,q}(z\0')} \,.
    $$
    This implies the function $f = \sum f^n$ satisfies 
    $$
    	\norm{f}_T \leq CM\norm{a\0}_{t_{\alpha}^{p,q}(z\0')} \sum_n (\delta CM)^n
    $$
    which is a convergent sum for $\delta < 1/[CM]$. Thus $f \in AT_{\alpha}^{p,q}$ as $AT_{\alpha}^{p,q} \subseteq\nobreak T_{\alpha}^{p,q}$ is closed by \thref{?16?}. Finally, we evaluate 
    $$
    	f(z_k') = \sum_{n=0}^{\infty} f^n(z_k') =  \sum_{n=1}^{\infty} a_k^{n} - a_k^{n+1} = a_k^0 = a_k
    $$
    as desired.
\end{proof}

We now come to the conclusion of this subsection.

\Proposition{?36?}{
 	If $z\0$ is interpolating for $T_{\alpha}^{p,q}$, then $D^+(z\0) < (1+\alpha)/q$.
}
\begin{proof}
	If $z\0$ is interpolating for $T_{\alpha}^{p,q}$ then so is every element of $W(z\0)$ by \thref{?34?}. Thus every element of $W(z\0)$ is a zero set for $T_{\alpha}^{p,q}$ and also for $A^{p}_{\beta + \epsilon}$ for $\beta = (1+\alpha)p/q - 1$ and $\epsilon>0$ by \thref{?20?}. Applying \thref{?32?} and the interpolation theorem for the Bergman space (\thref{?A?}), we obtain
    $$
    D^+(z\0) \leq ( [(1+\alpha)p/q - 1 + \epsilon] + 1 )/p = (1+\alpha)/q + \epsilon/p
    $$
    for each $\epsilon > 0$ giving
    $$
        D^+(z\0) \leq (1+\alpha)/q\,.
    $$
    
    Choose some sequences $r_n$ increasing to one and $w_n \in \dd$ so that
    $$
        D^+(z\0) - 1/n \leq D( z\0, w_n, r_n ) \,.
    $$
    Set $z\0^n = g_{w_n}^{-\eta}(z\0)$ as in \thref{?P?} where $\eta < \delta(z\0)/4$ is small so that \thref{?35?} implies $z\0^n$ is still interpolating as $\rho(z_k , z_k^n) < \eta$ for all $k$ by (1) of \thref{?P?}. From the above, we have $D^+(z\0^n) \leq (1+\alpha)/q$ for each $n$. By \thref{?P?}, there is $C>0$ depending only on $\delta$ such that we have
    \begin{align*}
    	D^+(z\0) &\leq D(z\0\SP,w_n,r) + 1/n \\&
    	\leq (1-\eta)\group{ D(z\0^n\,, w_n, r_n) - C[\log(1-r_n)]^{-1} } + 1/n \\&
    	\leq (1-\eta) \group{ (1+\alpha)/q - C[\log(1-r_n)]^{-1} } + 1/n
    \end{align*}
    which implies
    $$
    	D^+(z\0) \leq (1-\eta) (1+\alpha)/q < (1+\alpha)/q
    $$
    by letting $n \arrow \infty$ and noting that $(1-\eta) < 1$. This completes the proof. 
\end{proof}

Combining \thref{?26?,?36?} gives \thref{interpolationTheorem}, one of our main theorems.

\subsection{ Interpolation for Measures }

We now examine a related interpolation problem for measures. For a Borel measure $\mu$, we define
$$
    \scriptAT_{\alpha}^{p,q}(\mu) = \varinjlim AT_{\alpha}^{p,q}(U,\mu)
$$
where the limit is taken over all those open sets $U$ containing the support of $\mu$ and two functions are identified if they agree on $U \cap V$. If $\mu$ is a Carleson measure, then there is a natural map $AT_{\alpha}^{p,q} \arrow \scriptAT_{\alpha}^{p,q}(\mu)$ which is continuous. We write $[f]$ for the image of $f$ under this map. We will say that a measure $\mu$ is interpolating if this natural map is surjective and $\scriptAT_{\alpha}^{p,q} \neq AT_{\alpha}^{p,q}$. In the case that $\mu = \nu^{z\0}$ for $z\0$ uniformly discrete, we have that $t_{\alpha}^{p,q} = \scriptAT_{\alpha}^{p,q}(\mu) = T_{\alpha}^{p,q}(\mu)$ by representing any sequence $a\0$ as a function $F[a\0]$ as in \thref{?5?} where $$
    F[a\0] = \sum_{k=1}^{\infty} a_k \chi_{D(z_k,\epsilon)}
$$
for $\epsilon < \delta(z\0)$. Moreover, $F[a\0]$ is holomorphic on $U = \cup_k D(z_k,\epsilon)$. In this identification, the evaluation map is the natural map $AT_{\alpha}^{p,q} \arrow \scriptAT_{\alpha}^{p,q}(\nu^{z\0})$. Thus $\nu^{z\0}$ is interpolating for $AT_{\alpha}^{p,q}$ if and only if $z\0$ is an interpolating sequence for $AT_{\alpha}^{p,q}$. We now prove that the only interpolating Carleson measures $\mu$ for $AT_{\alpha}^{p,q}$ are those which are equivalent to $\nu^{z\0}$ for some uniformly discrete $z\0$ with $D^+(z\0) < (1+\alpha)/q$. (We say two Borel measures $\mu_1$ and $\mu_2$ are equivalent if there is $C$ such that 
$$
    C^{-1}\mu_2(E) \leq \mu_1(E) \leq C \mu_2(E)
$$
for each Borel set $E$.) For the remainder of the section, assume that $\mu$ is interpolating for $AT_{\alpha}^{p,q}$. 

\Lemma{?37?}{
    If $\mu$ is an interpolating Carleson measure, then the support of $\mu$ is a sequence.
}
\begin{proof}
    If the support of $\mu$ is all of $\dd$, we have that $\scriptAT_{\alpha}^{p,q}(\mu) = AT_{\alpha}^{p,q}(\dd, \mu)$. Hence, $AT_{\alpha}^{p,q} \arrow \scriptAT_{\alpha}^{p,q}(\mu)$ is bijective as $[f]_{\mu} = [g]_{\mu}$ implies that $f = g$ almost everywhere in $\mu$. As $f$ and $g$ are continuous, we must have $f=g$. This contradicts the assumption that $\mu$ is interpolating. By the open mapping theorem, the natural map is an isomorphism. Now assume that the support $E$ of $\mu$ misses a point. We show that the support $E$ of $\mu$ must be a set with no accumulation points in $\dd$. Assume that $E$ has an accumulation point $z_0$ in $\dd$. Let $z\0 \subset E$ converge to $z_0$. As $E$ is closed in $\dd$ but not all of $\dd$, there is $w_0 \in \dd$ so that $D(w_0,r) \subseteq \dd\excise E$. Define $f(z) = (z-w_0)^{-1}$. Then $f$ is holomorphic on $\dd \excise \set{w_0}$ and bounded on $E$. Thus $f$ defines some function $[f]$ in $\scriptAT_{\alpha}^{p,q}(\mu)$. Now assume $g$ is some analytic function on $\dd$ that projects to $[f]$. By definition, we must have $f=g$ on $E$. As $E$ has an accumulation point, $f=g$ on $\dd \excise \set{w_0}$. This is a contradiction. Hence $E$ can not have any accumulation points. This implies that $E$ must be a sequence. 
\end{proof}

Using the Carleson condition, we can force a bound on the measure of points. From here, we can use the interpolating condition to force the support of $\mu$ to be uniformly discrete.

\Lemma{?38?}{
    If $\mu$ is an interpolating Carleson measure for $AT_{\alpha}^p$ supported on $z\0$, then $z\0$ is uniformly discrete and $\mu(E) \leq C\nu(E)$ where $\nu = \nu^{z\0}$.
}
\begin{proof}
    Let $z\0$ be the support of $\mu$ and $r \in (0,1)$ be fixed. By the Carleson condition, 
    \[
        \mu(\set{z_k}) \leq \mu(D(z_k,r)) \leq C A( D(z_k,r) ) \leq C' (1-\abs{z_k}^2)^2 \tag{$\star$}
    \]
    where $C'$ depends only on $\mu$. Write $C(\mu)$ for this constant.
    We now obtain an interpolation constant for $\mu$ as in \thref{?30?}. Let $Z$ be the kernel of the natural map and let $I:AT_{\alpha}^{p,q}/Z \arrow \scriptAT_{\alpha}^{p,q}(\mu)$ be the quotient map. By definition of the quotient norm and the Carleson condition, $I$ is continuous. As $\mu$ is interpolating, I is bijective. By the open mapping theorem, the inverse of $I$ is continuous. This gives the interpolation constant $M = M(\mu)$. 
    
    We now follow the scheme of \thref{?31?} as follows: Define 
    $$
        Sf(z) = (1-\abs{z}^2)^{\gamma} \abs{f(z)}
    $$ for $\gamma = (1+\alpha)/q + 1/p$. For fixed $n$, let $a_k = \delta_k^n$ and $\nu = C(\mu)(1-\abs{z_n}^2)^2\delta_{z_n}$. As $a_k = 0$ for $k \neq n$, we can use ($\star$) above to get
    $$
        \norm{a\0}_{T_{\alpha}^{p,q}(\mu)} \leq \norm{a\0}_{T_{\alpha}^{p,q}(\nu)} \leq C (1-\abs{z_n}^2)^{\gamma}
    $$
    by direct calculation. Find $f \in AT_{\alpha}^{p,q}$ such that 
    $$
        \norm{f}_{T_{\alpha}^{p,q}} \leq M \norm{a\0}_{T_{\alpha}^{p,q}(\mu)} \leq CM(1-\abs{z_k}^2)^{\gamma} \,.
    $$
    Assume that $\rho(z_k,z_n)$ is small enough and apply \thref{?29?} together with a direct calculation to obtain
    $$
        (1-\abs{z_n}^2)^{\gamma} = \abs{Sf(z_n)-Sf(z_k)} \leq \rho(z_k,z_n) CM(1-\abs{z_n}^2)^{\gamma}
    $$
    giving a lower bound on $\rho(z_k,z_n)$. If \thref{?29?} does not apply, then we know $\rho(z_k,z_n) \geq R > 0$ for some fixed $R$. Thus $z\0$ is uniformly discrete and ($\star$) implies that $\mu(E) \leq C(\mu) \nu(E)$ as desired.
\end{proof}

The above lemma implies that $\scriptAT_{\alpha}^{p,q}(\nu^{z\0})$ includes into $\scriptAT_{\alpha}^{p,q}(\mu)$ where $z\0$ is the support of $\mu$. We now show it is also surjective.

\Proposition{?39?}{
    If $\mu$ is an interpolating Carleson measure, then $\scriptAT_{\alpha}^{p,q}(\mu) = t_{\alpha}^{p,q}(z\0)$ for some uniformly discrete sequence.
}

\begin{proof}
    As we noted, the assumptions imply that $t_{\alpha}^{p,q}(z\0)$ includes into $\scriptAT_{\alpha}^{p,q}(\mu)$ where $z\0$ is the support of $\mu$. However, the natural map $AT_{\alpha}^{p,q} \arrow \scriptAT_{\alpha}^{p,q}(\mu)$ must factor through $t_{\alpha}^{p,q}(z\0) = \scriptAT_{\alpha}^{p,q}(\nu^{z\0})$ by \thref{?14?}. Thus if there were some sequence $a\0 \in \scriptAT_{\alpha}^{p,q}(\mu)$ which is not an element of $t_{\alpha}^{p,q}(z\0)$, then $a\0$ can not be in the image of the natural map. This contradicts the assumption that $\mu$ is interpolating.
\end{proof}

We now collect everything into one theorem:

\Theorem{?40?}{
    A Carleson measure $\mu$ is interpolating for $AT_{\alpha}^{p,q}$ if and only if $\mu$ is equivalent to $\nu^{z\0}$ for some uniformly discrete sequence $z\0$ with $D^+(z\0) < (1+\alpha)/q$. 
}
\begin{proof}
    We have already seen that $\nu^{z\0}$ is interpolating if $z\0$ is an interpolating sequence. As equivalent measures yield isomorphic tent spaces, we have that $\mu$ equivalent to $\nu^{z\0}$ implies that $\mu$ is interpolating. 
    
    Now assume that $\mu$ is an interpolating Carleson measure. Then \thref{?39?} implies that $t_{\alpha}^{p,q} = \scriptAT_{\alpha}^{p,q}(\mu)$. Picking $a_k = \delta_k^n$, we calculate
    $$
        \norm{a\0}_{T_{\alpha}^{p,q}(\mu)} = \mu(\set{z_k})^{1/q} (1-\abs{z_k}^1)^{(\alpha-1)/q+1/p}
    $$
    and
    $$
        \norm{a\0}_{T_{\alpha}^{p,q}(\mu)} = (1-\abs{z_k}^1)^{(\alpha+1)/q+1/p}
    $$
    so that $\norm{a_k}_{T_{\alpha}^{p,q}(\mu)} \simeq \norm{a\0}_{t_{\alpha}^{p,q}}$ implies that $\mu$ is equivalent to $\nu^{z\0}$. It is clear that $z\0$ is interpolating for $AT_{\alpha}^{p,q}$ so that $D^+(z\0) < (1+\alpha)/q$ as desired.
\end{proof}

The above result justifies why we have chosen to consider $t_{\alpha}^{p,q}$. It may be interesting to ask what happens for interpolating measures which are not Carleson (if such measures exist). We do not see a clear way to get a handle on this now.

\section{Sampling}
We now turn to the sampling problem. To show that the evaluation operator $E$ is injective with closed range, we prove that $E$ is bounded below. We prove there is a constant $C$ such that the inequalities
$$
    C^{-1} \norm{f}_{T_{\alpha}^{p,q}} \leq \norm{f}_{T_{\alpha}^{p,q}(\nu^{z\0})} \leq C \norm{f}_{T_{\alpha}^{p,q}}
$$
hold for every $f \in T_{\alpha}^{p,q}$. For convenience, we recall the main theorem of this section. 

\MainTheorem{samplingTheorem'}{
    A sequence $z\0 \subseteq \dd$ is sampling for $AT_{\alpha}^{p,q}$ if and only if $z\0$ is a Carleson sequence and there is a uniformly discrete subsequence $z\0'$ of $z\0$ such that $D^-(z\0') > (1+\alpha)/q$.
}
\subsection{Sufficiency of the Sampling Condition}

We first define constants that quantify \enquote{how well} a sequence samples for the  the tent spaces. There are analogous constants for sampling in Bergman spaces.

\Definition{sampling_constant}{
    Suppose $z\0$ is sampling for $T_{\alpha}^{p,q}$ such that $L$ and $U$ are the best constants with
    $$
        L \norm{f}_{T_{\alpha}^{p,q}} \leq \norm{f}_{T_{\alpha}^{p,q}(\nu^{z\0})} \leq U \norm{f}_{T_{\alpha}^{p,q}}
    $$
    for all $f \in AT_{\alpha}^{p,q}$. We then call $L = L(z\0)$ the \emph{lower sampling constant} and $U = U(z\0)$ the \emph{upper sampling constant}. We can similarly define the upper and lower sampling constants for the Bergman space by replacing the tent space with the Lebesgue space in the above. We similarly define $\tilde{L}(z\0)$ and $\tilde{U}(z\0)$ by replacing $T_{\alpha}^{p,q}$ with $\tilde{T}_{\alpha}^{p,q}$ in the above inequalities.
}

Applying \thref{?13?} to the measure $\nu^{z\0}$, we know that the inequality
$$
    \norm{ f(z\0) }_{t_{\alpha}^{p,q}} \leq U \norm{f}_{T_{\alpha}^{p,q}}
$$
always holds under our assumption that $z\0$ is Carleson. In order to prove the other sampling inequality, we hope to estabilish inequalities of the form
$$
    \abs{f(z)}^s \leq S_{b,c}[\gap{\abs{f(z\0)}^s}] (z)
$$
where $S_{b,c}$ is the key operator from \thref{?25?}. Note that the sampling condition implies that $z\0$ is sampling for some related Bergman space. We can use this to get an inequality as above. However, we first need the following lemma.

\Lemma{?41?}{
    Let $w \in \dd$. If $z\0$ is sampling for the Bergman space $A_{\alpha}^p$ then $\phi_w(z\0)$ is still sampling for $A_{\alpha}^p$ with the same sampling constants.
}
The proof of the above is well understood and follows from the methods of \cite{DurenSchuster}. We also include this related lemma which we use later. Its proof is virtually identical to the above.

\Lemma{?42?}{
    Let $w \in \dd$. If $z\0$ is sampling for the tent space $AT_{\alpha}^{p,q}$ then $\phi_w(z\0)$ is still sampling for $AT_{\alpha}^{p,q}$. Moreover, $L(z\0) \simeq  L(\phi_w(z\0))$ and $U(z\0) \simeq U(\phi_w(z\0))$ with equivalence constant independent of $w$.
}
\begin{proof}
    Define $S_wf(z) = f\comp\phi_w(z) \dott \phi_w'(z)^{\gamma}$ for $\gamma = 1/p + (1+\alpha)/q$ so that $S_w^2 = \id$. By \thref{?10?,?11?} we have 
    $$
        \norm{S_wf}_{T_{\alpha}^{p,q}} \simeq \norm{f}_{T_{\alpha}^p}
        \andEq
        \norm{S_wf}_{T_{\alpha}^{p,q}(\nu)} \simeq \norm{f}_{T_{\alpha}^{p,q}\group{\nu_w}}
    $$
    where $\nu = \nu^{z\0}$ and $\nu_w = \nu^{\phi_w(z\0)}$. Applying the sampling criterion to $S_wf(z)$ we get 
    $$
        L(z\0) \norm{S_wf(z)}_{T_{\alpha}^{p,q}} \leq \norm{S_wf}_{T_{\alpha}^{p,q}\group{\nu}} \leq U(z\0) \norm{S_wf}_{T_{\alpha}^{p,q}} 
    $$
    which becomes
    $$
        L(z\0) \norm{f(z)}_{T_{\alpha}^{p,q}} \lesssim \norm{f}_{T_{\alpha}^{p,q}\group{\nu_w} } \lesssim U(z\0) \norm{f}_{T_{\alpha}^{p,q}}
    $$
    upon application of \thref{?10?,?11?}. This holds for all $f \in T_{\alpha}^{p,q}$ so that $L(z\0) \leq C_1 L(\phi_w(z\0))$ and $U(\phi_w(z\0)) \leq C_2 U(z\0)$. Thus $\phi_w(z\0)$ is sampling. Applying the same reasoning to the sampling sequence $\phi_w(z\0)$ gives $L(\phi_w(z\0)) \leq C_3 L(z\0)$ and $U(z\0) \leq C_4 U(\phi_w(z\0))$ so that $L(z\0) \simeq L(\phi_w(z\0))$ and $U(z\0) \simeq U(\phi_w(z\0))$ with equivalence constant independent of $w$.
\end{proof}

We now apply the above results to obtain the following estimates on the size of $f$. Note that it is critical the sampling constants of $\phi_{\zeta}(z\0)$ do not depend on $\zeta$. (The author wants to thank Daniel Luecking for motivating the following observation.)

\Proposition{?43?}{
    Suppose $z\0$ is sampling for $A_{c}^{s}$. We have the estimate
    $$
        \abs{f(z)}^{s} \leq C (1-\abs{z}^2)^{b} \sum_{k=1}^{\infty} \abs{f(z_k)}^{s} \frac{ (1-\abs{z_k}^2)^{c+2}}{\abs{1-\overbar{z}z_k}^{c +b + 2} } = S_{b,c}[\gap{\abs{f(z_k)}^s}](z)
    $$
    for any real $b$. Here, $S$ is taken from \thref{?25?}.
}

\begin{proof}
    Set $w_k = \phi_z(z_k)$ so that $w\0$ is still sampling with the same sampling constants by \thref{?41?}. The lower sampling inequality and boundedness of point evaluations give
    $$
        \abs{h(0)}^{s} \leq \norm{h}_{A_{c}^{s}} \leq C \sum_{k=1}^{\infty} \abs{h(w_k)}^{s} (1-\abs{w_k}^2)^{c+2} \,.
    $$
    Note that $C$ depends on the sampling constant of $w\0\SP$. Apply this to $h = g\comp\phi_z$ to get
    $$
        \abs{g(z)}^{s}  \leq 
        C \sum_{k=1}^{\infty} \abs{g(z_k)}^{s} (1-\abs{\phi_z(z_k)}^2)^{c+2} 
        = C (1-\abs{z}^2)^{2+c} \sum_{k=1}^{\infty} \abs{g(z_k)}^{s} \frac{ (1-\abs{z_k}^2)^{c+2}}{\abs{1-\overbar{z}z_k}^{2(c + 2)} }
    $$
    by standard formulas for $\phi_z(w)$. Finally, we apply the above to the analytic function
    $$
        g_z(w) = f(w) \groupFrac{1-\abs{z}^2}{1-\overbar{z}w}^{(b-c-2)/s}.
    $$
    Noting that $g_z(z) = f(z)$, we get the result.
\end{proof}

We now leverage the above result together with the boundedness result, \thref{?25?}, to obtain sampling --- conveniently, in our applications, $c+1$ will be greater than $(1+\alpha)/q$. Note that we must apply some of the same tricks as in \thref{?10?} to deal with the case when $p,q<1$. Let us define a supersequence $z\0$ of $z\0'$ to be a sequence such that $z\0'$ is a subsequence of $z\0\SP$.

\Proposition{?44?}{
    If $z\0$ is a Carleson sequence and has a uniformly discrete subsequence $z\0'$ such that $D^-(z\0')>(1+\alpha)/q$, then $z\0$ is sampling for $AT_{\alpha}^{p,q}$.
}
\begin{proof}
    As a Carleson supersequence of a sampling sequence is still sampling, it suffices to take $z\0 = z\0'$\,. Moreover as $z\0$ is Carleson,  the mapping $f \mapsto f(z\0)$ is continuous from $T_{\alpha}^{p,q} \arrow t_{\alpha}^{p,q}(z\0)$ by \thref{?13?}. Thus it suffices to show the lower sampling inequality
    $$
        \norm{f}_{T_{\alpha}^{p,q}} \leq C \norm{f(z\0)}_{t_{\alpha}^{p,q}}
    $$
    on a dense subset of $AT_{\alpha}^{p,q}$, namely polynomials. (See corollary 8 of \citePage{18}{Perala} for the proof that polynomials are dense in $AT_{\alpha}^{p,q}$.)
    Define $S$ as in \thref{?25?} by
    $$
        S_{b,c}f(z) = (1-\abs{z}^2)^{b} \sum_{k=1}^{\infty} \abs{f(z_k)} \frac{ (1-\abs{z_k}^2)^{c+2}  }{\abs{1-\overbar{z_k}z}^{c+b+2} }
    $$ 
    Let $s < \Min{p,q}$ so that $P=p/s$ and $Q=q/s$ are bigger than $1$. Define 
    $$
        Rf(z) = \abs{f(z)}^{s}
    $$
    and choose $c$ so that $(1+c)/s \in ((1+\alpha)/q, D^-(z\0) )$. For $f$ bounded, \thref{?43?} implies 
    $$
        Rf(z) \leq C S_{b,c}[Rf(z\0)](z)
    $$
    for any $b \in \rr$ as $z\0$ is sampling for $A_{c}^s$. Moreover, direct calculation gives that 
    $$
        \norm{f}_{T_{\alpha}^{p,q}(\mu)}^p = \norm{Rf}_{T_{\alpha}^{P,Q}(\mu)}^P
    $$
    for each measure $\mu$. As $c$ satisfies $(1+c)/s > (1+\alpha)/q$, we know 
    $$
        c +1 > (1+\alpha)/(q/s) = (1+\alpha)/Q.
    $$
    Thus by \thref{?25?} we have $S_{b,c} : t_{\alpha}^{P,Q}(z\0) \arrow T_{\alpha}^{P,Q}$ is bounded for $b$ sufficiently negative. Recall that $t_{\alpha}^{p,q}(z\0) = T_{\alpha}^{p,q}(\nu)$ for $\nu = \nu^{z\0}$. Thus we get
    $$
        \norm{f}_{T_{\alpha}^{p,q}}^p = \norm{Rf}_{T_{\alpha}^{P,Q}}^P \leq C \norm{S_{b,c}Rf(z\0)}_{T_{\alpha}^{P,Q}}^P \leq C' \norm{Rf}_{T_{\alpha}^{P,Q}\group{\nu}}^P  = C' \norm{f}_{T_{\alpha}^{p,q}\group{\nu}}^p
    $$
    giving the sampling result.
\end{proof}

Note that we needed $P,Q > 1$. In general, the author does not believe that the operator $S:t_{\alpha}^{P,Q} \arrow T_{\alpha}^{P,Q}$ as above is bounded when $P,Q \leq 1$. Note that we have additional restrictions on the parameters of $S$ when $P,Q \leq 1$. (Refer to \thref{?23?,?24?} for more information.)

\subsection{Necessity of the Sampling Condtion}

As in the interpolation problem, we will need to perturb a sampling sequence to obtain a strict inequality. Without perturbation, we can only prove that
$$
    D^-(z\0') \geq (1+\alpha)/q
$$
for some subsequence $z\0'$ of the sampling sequence, $z\0\SP$. In order to state the precise result on perturbation, we need the following definitions. We also use this technology to extract a uniformly discrete subsequence $z\0'$ of a sampling sequence with $z\0'$ is still being sampling.

\Definition{carleson_class}{
    Given $\delta>0$ and $N\in\nn$, we define the \emph{Carleson class}, $\scriptD(\delta,N)$, to be the collection of all sequences which can be written as a union of $N$ uniformly discrete sequences which all have separation constants bounded below by $\delta$.
}

Inside of every Carleson class, we have a sampling class as follows:

\Definition{sampling_class}{
    Given $\delta>0$, $L>0$, and $N \in \nn$, we define the \emph{sampling class}, $\scriptS(\delta,N,L)$, to be those $z\0 \in \scriptD(\delta,N)$ such that $L(z\0) \geq L$.
}

The goal is to show that two sequences $z\0$ and $w\0$ which live in the same Carleson class $\scriptD(\delta,N)$ satisfy the following general principle: Whenever $z\0$ and $w\0$ are \enquote{close} to one another, then $z\0$ and $w\0$ live in \emph{neighboring} sampling classes. The notion of \enquote{closeness} is given by the Hausdorff distance:

\Definition{hausdorff_distance}{
    The \emph{pseudohyperbolic Hausdorff distance}, $[A,B]$, for closed sets $A$ and $B$ is defined by
    $$
        [A,B] = \Inf{ \delta > 0 : B \subseteq N_{\delta}(A) \AND A \subseteq N_{\delta}(B) }
    $$
    where $N_{\delta}$ is the pseudohyperbolic neighborhood.
}
Before we continue, we need a result about the Bochner-Lesbegue $p$-space,  
$$
    \scriptL^{p,q} = L^p(\bd\dd  \myspace ; \myspace \ell^q) = \set{ f: \bd\dd \arrow \ell^q \myspace[2]\mathlarger{\mid}\myspace[2] \norm{f(\dott)}_{\ell^q} \in L^p },
$$
of functions from $\bd\dd$ to $\ell^q$. The norm is given by the formula
$$
    \norm{f}_{\scriptL^{p,q}}^p = \int_{\bd\dd} \norm{ f(\zeta) }_{\ell^q}^p \dif l(\zeta) \,.
$$
The following will allow us to estimate 
$$
    \bigabs{ \norm{f}_{t_{\alpha}^{p,q}(z\0)} - \norm{f}_{t_{\alpha}^{p,q}(w\0)} }
$$
under assumptions about $z\0$ and $w\0\SP$:
\Lemma{?45?}{
     If $s = \Min{1,p,q}$, then the function $\norm{\dott} \defeq \norm{\dott}_{\scriptL^{p,q}}^s$ is an $s$-norm. 
}
\noindent The proof of the above is the same as \thref{?4?}. We omit the details. We now estimate
$$
    \bigabs{ \norm{f}_{t_{\alpha}^{p,q}(z\0)} - \norm{f}_{t_{\alpha}^{p,q}(w\0)} }
$$
whenever $w\0$ is a small perturbation of $z\0$.
\Lemma{?46?}{ 
    Fix $\delta>0$, $N \in \nn$, and $\lambda > \Max{1,q/p}$. Define the spaces 
    $$
        Z=\myTilde{T}_{\alpha,\lambda}^{p,q}\group{\nu^{z\0}} 
        \andEq 
        W=\myTilde{T}_{\alpha,\lambda}^{p,q}(\nu^{w\0}).
    $$ 
    If $z\0\SP,w\0 \in \scriptD(\delta,N)$ with 
    $$
        \sup_k\set{\rho(z_k, w_k)} = \Delta < \delta/8,
    $$
    then there is a positive function $M(x)$, depending only on $N$, with $M(x) \arrow 0$ monotonically as $x \arrow 0$ satisfying
    $$
        \bigabs{\norm{f}_Z^s - \norm{f}_W^s } < M(\Delta) \norm{f}_{T_{\alpha}^{p,q}}^s
    $$
    where $s = \Min{p,q,1}$.
}
\begin{proof} 
    In the subsequent, $C$ will represent some positive constant depending only on $N$ and $\delta$.
    Write $Sf(z) = (1-\abs{z}^2)^{(1+\alpha)/q}\abs{f(z)}$. We also define the $\ell^p$-valued functions $a$, $b$, and $c$ by
    \begin{align*}
        a_k(\zeta) &= K(z_k,\zeta)^{\lambda/q} Sf(z_k) \\ 
        b_k(\zeta) &= K(z_k,\zeta)^{\lambda/q} Sf(w_k) \\
        c_k(\zeta) &= K(w_k,\zeta)^{\lambda/q} Sf(w_k)
    \end{align*} 
    for 
    $$
        K(z,\zeta) = \frac{1-\abs{z}^2}{\abs{1-\overbar{\zeta}z}}
    $$
    as in \thref{?C?}. By definition, we have the following identities:
    \[
        \norm{f}_Z = \norm{ a }_{\scriptL^{p,q}}
        \andEq
        \norm{f}_W = \norm{ c }_{\scriptL^{p,q}} \tag{$\star$} \,.
    \]
    Let $\norm{\Dott}_p = \norm{\Dott}_{L^p(\bd\dd)}$. For $D_k = D(z_k,\Delta/2)$, we apply \thref{?28?} to get
    $$
        \abs{ a_k(\zeta) - b_k(\zeta) }^q  \leq
        \Delta C \hspace{1pt} K(z_k,\zeta)^{\lambda}  \int_{D_k} \abs{f(z)}^q \dif\sigma_{\alpha}(z) \leq
        \Delta C' \int_{D_k} K(z,\zeta)^\lambda \abs{f(z)}^q \dif\sigma_{\alpha}(z)
    $$
    using \thref{?G?} for the last inequality. We then use \thref{?45?} with the above to get
    $$
        \bigabs{ \norm{a}_{\scriptL^{p,q}}^s - \norm{b}_{\scriptL^{p,q}}^s } \leq 
        \norm{ a - b }_{\scriptL^{p,q}}^s \leq
        C^s\Delta^{s}\norm*{ \sum_k \int_{D_k} K(z,\Dott)^{\lambda} \abs{f(z)}^q
        \dif\sigma_{\alpha}(z) }_{p/q}^{s/q} \leq
        C^s\Delta^s\norm{f}_{T_{\alpha}^{p,q}}\,.
    $$
    The last inequality above follows as $\Delta < \delta$ implies the sets $\set{D_k}$ are disjoint. As in page 42 of \cite{DurenSchuster} and page 87 of \cite{Zhu}, we have the estimate
    $$
        \abs*{1-\frac{K(z_k,\zeta)^{\lambda}}{K(w_k,\zeta)^{\lambda}} } < M_1(\Delta)
    $$ 
    for each $k$ and $\zeta$  where $M_1(x)$ depends only on $N$ and $M_1(x) \arrow 0$ monotonically as $x \arrow 0$. By definition of $b$ and $c$ together with the above estimate, we have that
    $$
        \abs{c_k(\zeta)-b_k(\zeta)}^q = 
        \abs{b_k(\zeta)}^q{ \abs*{1-\frac{K(w_k,\zeta)^{\lambda}}{K(z_k,\zeta)^{\lambda}} } }^q \leq
        M_1(\Delta) K(z_k,\zeta)^{\lambda} (1-\abs{z_k}^2)^{1+\alpha}\abs{f(z_k)}^q \,.
    $$
    From this, \thref{?45?}, and estimates from \thref{?G?,?J?}; we see that
    $$
        \bigabs{ \norm{b}_{\scriptL^{p,q}}^s - \norm{c}_{\scriptL^{p,q}}^s } \leq 
        \norm{ b - c }_{\scriptL^{p,q}}^s \leq
        C M_1(\Delta)^s\norm*{ \sum_k \int_{D_k} K(z,\Dott)^{\lambda} \abs{f(z)}^q
        \dif\sigma_{\alpha}(z) }_{p/q}^{s/q} \,.
    $$
    As the sets $\set{D_k}$ are disjoint, the above is bounded by $CM_1(\Delta)^s\norm{f}_{T_{\alpha}^{p,q}}^s$\,. Summarizing, we have the two estimates
    $$
        \bigabs{ \norm{a}_{\scriptL^{p,q}}^s - \norm{b}_{\scriptL^{p,q}}^s } \leq C \Delta^s \norm{f}_{T_{\alpha}^{p,q}} \andEq \bigabs{ \norm{b}_{\scriptL^{p,q}}^s - \norm{c}_{\scriptL^{p,q}}^s } \leq C M_2(\Delta) \norm{f}_{T_{\alpha}^{p,q}}
    $$
    By the triangle inequality and the identities $(\star)$, we have
    $$
        \bigabs{ \norm{f}_Z^s - \norm{f}_W^s } \leq
        \bigabs{ \norm{a}_{\scriptL^{p,q}}^s - \norm{b}_{\scriptL^{p,q}}^s } + \bigabs{ \norm{b}_{\scriptL^{p,q}}^s - \norm{c}_{\scriptL^{p,q}}^s } \leq
        M(\Delta) \norm{f}_{T_{\alpha}^{p,q}}^s
    $$
    with $M(x) = C\bracket[\big]{ x^s + M_1(x)^s }$ going to zero monotonically as $x \arrow 0^+$.
\end{proof}

We now prove our \enquote{general principle} in two steps. We start by assuming that $w\0$ is a perturbation of $z\0\SP$. We then generalize to the case when $[z\0\SP,w\0]$ is sufficiently small.

\Lemma{?47?}{
    Let $L$ and $\delta$ be positive constants. For each natural number $N$, there are constants $R \in (0, \delta/8)$ and $C>0$ satisfying the following: Whenever $z\0 \in \scriptS(\delta,N,L)$ and $w\0 \in \scriptD(\delta,N)$ with 
    $$
        \Sup{\rho(z_k,w_k)} = \Delta < R
    $$ then $L(w\0) \geq C L(z\0)$ for $C$ not depending on $z\0$ or $w\0$.
}
\begin{proof}
    For $\lambda > \Max{1,q/p}$, we will compute the lower sampling constant with respect to the equivalent norm of $\myTilde{T}_{\alpha,\lambda}^{p,q}$ as the calculations are easier. $M(x)$ will represent a function such that $M(x) \arrow 0$ monotonically as $x \arrow 0$. $M$ is allowed to depend on $\delta$, $N$, and $L$. Set $Z=\myTilde{T}_{\alpha,\lambda}^{p,q}(\nu^{z\0})$ and  $W=\myTilde{T}_{\alpha,\lambda}^{p,q}(\nu^{w\0})$. Recall that
    $$ 
        \myTilde{L}(z\0) = \Sup[\Big]{C : C\norm{f}_{\myTilde{T}_{\alpha,\lambda}^{p,q}} \leq \norm{f}_{\myTilde{T}_{\alpha,\lambda}^{p,q}\group[\gLarger]{\nu^{z\0}}} , \forall f \in \myTilde{T}_{\alpha,\lambda}^{p,q} }
    $$
    is the lower sampling constant with respect to $\myTilde{T}_{\alpha,\lambda}^{p,q}$. Note that $L(\dott) \simeq \myTilde{L}(\dott)$ by \thref{?C?}. For $\epsilon>0$ and any $f \in \myTilde{T}_{\alpha,\lambda}^{p,q}$ of norm 1, we have 
    $\norm{f}_Z \geq \myTilde{L}(z\0)$. Moreover, we can choose $f$ of norm 1 so that $\myTilde{L}(w\0) \geq \norm{f}_W-\epsilon$ for $\epsilon>0$. For $s = \Min{p,q,1}$, the above gives that
    $$
        \norm{f}_W^s \geq \norm{f}_Z^s - \bigabs{ \norm{f}_Z^s - \norm{f}_W^s } \geq \group{\myTilde{L}(z\0)-\epsilon}^s - M(\Delta)^s
    $$
    where $M$ is from \thref{?46?}. Taking $R$ small enough depending on $N$ and $\delta$, we have
    $$
        M(\Delta) < \myTilde{L}(z\0)^s/2
    $$
    as $\Delta < R$. Recalling that $\myTilde{L}(w\0)\geq \norm{f}_W-\epsilon$ and letting $\epsilon \arrow 0$  gives 
    $$
        \myTilde{L}(w\0)^s \geq \myTilde{L}(z\0)^s - M(\Delta) \geq \myTilde{L}(z\0)^s/2
    $$ 
    so that $L(w\0) \geq C L(z\0)$ as desired.
    
\end{proof}

\noindent We are able to improve the above result to sequences which are close in the Hausdorff distance. (The author wants to thank Daniel Luecking again for the improved exposition in the following.)

\Proposition{?48?}{
    Let $[\dott,\dott]$ be the Hausdorff distance and $L>0$, $\delta>0$ and $N \in \nn$ be fixed constants. There are a constants $R \in (0, \delta/8)$ and $C>0$ such that whenever $z\0 \in \scriptS(\delta,N,L)$ and $w\0 \in \scriptD(\delta,N)$ with $[z\0\SP,w\0] = \Delta < R$ we have that $L(w\0) \geq C L(z\0)$ where $C$ is independent of $z\0$ and $w\0$.
}
\begin{proof} %See draft 8.1 for the old proof which was wrong
    We want to construct a sequence $w\0'$ satisfying $\norm{\dott}_{t_{\alpha}^{p,q}(w\0')} \leq N \norm{\dott}_{t_{\alpha}^{p,q}(w\0)}$ and $\rho(z_k,w_k') < \Delta$. We construct $w\0'$ as follows: Let $\hat{w}\0$ be the sequence $w\0$ with each point repeated $N$-times. 
    As $[z\0\SP,\hat{w}\0] = [z\0\SP,w\0]$, there is an injective function $k \mapsto n_k$ so that $\rho(z_k,\hat{w}_{n_k}) \leq \Delta$. Define $w_k' = \hat{w}_{n_k}$. As $z\0$ and $w\0'$ satisfy the hypothesis of \thref{?47?}, we see $L(w\0') \geq CL(z\0)$. By construction,
    $$
        N \norm{f}_{t_{\alpha}^{p,q}(w\0)} = \norm{f}_{t_{\alpha}^{p,q}(\hat{w}\0)} \geq \norm{f}_{t_{\alpha}^{p,q}(w\0')}. 
    $$
    Hence, $L(w\0) \geq (1/N)L(w\0') \geq (C/N) L(z\0)$ as desired.
    
\end{proof}
The above \enquote{general principle} is analogus to the result of \cite{DurenSchuster} on page 198. We now use the above to extract a uniformly discrete subsequence $z\0'$ of a sampling sequence $z\0$ so that $z\0'$ is still sampling.
\Corollary{?49?}{
    Every sampling sequence has a subsequence which is uniformly discrete and sampling.
}
\begin{proof}
    By \thref{?48?}, we choose some $R$ such that $[z\0\SP,w\0]<R$ implies $w\0$ is also sampling whenever $w\0$ is a uniformly discrete sequence. Now define $w\0$ inductively as follows: Take $w_1 = z_1$. After choosing the first $n$ elements of $w\0$ we define $w_{n+1}$ to be the first element of $z\0$ following $w_n$ such that $\rho(w_k,w_{n+1}) > R/4$ for each $k \leq n$. Then for each $z_k$ there must be some $l$ such that $\rho(z_k,w_l) < R$. To see this, assume the negation. Then $D(z_k,R/4) \cap w\0 = \emptyset$ but $z_k$ was never chosen. This is a contradiction. Thus $[z\0\SP,w\0] < R$ and $w\0$ is discrete by construction. Thus we are done.
\end{proof}

We need the following to estimate $D^-(z\0)$:

\Lemma{?50?}{
    Define $\Gamma_{\zeta}^r = \Gamma_{\zeta} \excise D(0,r)$. For a uniformly discrete sequence $z\0\SP$, we have the following estimate
    $$
        \sum_{z_k \in \Gamma_{\zeta}^r} (1-\abs{z_k}^2)^{1+\alpha} \leq C(1-r)^{1+\alpha}
    $$
    as $r \arrow 1$ for each $\zeta \in \bd\dd$.
}
\begin{proof}
    Take $r>3/4$. As $z\0$ is uniformly discrete, we may take $\delta \leq \delta(z\0)/8$ so that $\set{D(z_k,\delta)}$ are disjoint. By \thref{?6?}, we can choose a Stoltz family $\Lambda_{\zeta}$ large so that $z_k \in \Gamma_{\zeta}^r$ implies $D(z_k,\delta) \subseteq \Lambda_{\zeta}^R$ where $R = (r-\delta)/(1-r\delta)$. (Note that $R$ has been chosen with respect to the strong triangle inequality on the pseudohyperbolic metric.) We estimate 
    $$
        \sum_{z_k \in \Gamma_{\zeta}^r} (1-\abs{z_k}^2)^{1+\alpha} \leq C \int_{\Lambda_{\zeta}^R} \dif\sigma_{\alpha}
    $$
    by \thref{?G?}. We choose $r$ large enough and $\delta$ small enough that $R$ is larger than the aperture of $\Lambda$. There is a function, $\Theta(s)$, so that 
    $$
        \Lambda_{\zeta}^R = \zeta \dott \hspace{2pt} \set{ se^{i\theta} : s \in (R,1) \AND \abs{\theta} \in [0,\Theta(s) ) } \,.
    $$
    We now compute that
    $$
        \int_{\Lambda_{\zeta}^R} \dif\sigma_{\alpha} = C_{\alpha}\int_{R}^1 \int_{-\Theta(s)}^{\Theta(s)} s (1-s^2)^{\alpha-1} \dif \theta \dif s = 2C_{\alpha}\int_{R}^1 s (1-s^2)^{\alpha-1} 2\pi \dott l(I_s) \dif s
    $$
    where $I_s = \set{ \zeta \in \bd \dd : s \in \Lambda_{\zeta} } $ is as in \thref{?1?}. We conclude that
    $$
         \int_{\Lambda_{\zeta}^R} \dif\sigma_{\alpha}  = 2C_{\alpha}\int_{R}^1 s (1-s^2)^{\alpha-1} 2\pi \dott l(I_s) \dif s \leq C \int_R^1 s(1-s^2)^{\alpha} \dif s \leq C' (1-R^2)^{1+\alpha}
    $$
    by \thref{?1?} followed by a direct computation. Finally, $\delta < 1/8$ implies that $1-R^2 \leq C( 1-r^2)$ which gives our result as $1-r \simeq 1-r^2$ for $r \in [0,1]$.
\end{proof}

We now deduce the sampling result.

\Proposition{?51?}{
    If $z\0$ is sampling for $T_{\alpha}^{p,q}$ then $z\0$ is a Carleson sequence and there is a uniformly discrete subsequence $w\0$ of $z\0$ such that $D^-(w\0) > (1+\alpha)/q$.
}
\begin{proof}
    We first note that we are following the ideas of pages 205-207 in \cite{DurenSchuster}. By the upper sampling inequality, the inclusion 
    $$
        AT_{\alpha}^{p,q} \arrow AT_{\alpha}^{p,q}(\nu^{z\0}) \subseteq t_{\alpha}^{p,q}(z\0)
    $$
    is bounded. Thus, \thref{?13?} implies that every sampling sequence is a Carleson sequence. By the Carleson condition (\thref{?G?}) and \thref{?49?}, we can assume that $z\0$ is a uniformly discrete sampling sequence with $\delta = \delta(z\0)$. We only need to show $D^-(z\0) > (1+\alpha)/q$. Take $\epsilon_j \arrow 0$ with $\epsilon_j \in (0,1/2)$ and let $r_j \geq 1-\epsilon_j$. Choose $\zeta_j \in \dd \excise z\0$ such that $D(z\0^j,0,r_j) < D^-(z\0)+\epsilon_j$ where $z\0^j = \phi_{\zeta_j}(z\0)$. Note that $\delta(z\0^j) = \delta(z\0)$ for each $j$. By \thref{?42?}, $z\0^j$ is sampling with lower sampling constant bounded in terms of $L(z\0)$. Define the function $f^j$ by the finite product
    $$
        f^j(z)=\prod_{ k \in K_j  } \frac{\phi_{z_k^j}(z)}{\abs{z_k^j}} 
    $$
    where $K_j = \set{k \,:\,0 < \abs{z_k^j} < r_j}$ so that $f^j(z_k^j) = 0$ when $\abs{z_k^j} < r_j$. We have
    $$
        1 = \abs{f^j(0)} \leq C \norm{f^j}_{T_{\alpha}^{p,q}}
    $$
    by \thref{?16?} where $C$ is independent of $j$. Defining $\Gamma_{\zeta}^r = \Gamma_{\zeta}\excise D(0,r)$, we estimate that
    \begin{align*}
        \sum_{z_k^j \in \Gamma_{\zeta}} (1-\abs{z_k^j}^2)^{1+\alpha}\abs{f^j(z_k^j)}^q &=
        \sum_{z_k^j \in \Gamma_{\zeta}^r} (1-\abs{z_k^j}^2)^{1+\alpha}\abs{f_j(z_k^j)}^q \\&\leq
        \sum_{z_k^j \in \Gamma_{\zeta}^r} (1-\abs{z_k^j}^2)^{1+\alpha} \prod_{ 0<\abs{z_l^j} < r_j } \frac{1}{\abs{z_l^j}^q} \\&\leq
        C \group[\Bigg]{\myspace \prod_{ 1/2 < \abs{z_l^j} < r_j } \frac{1}{\abs{z_l^j}^q} } \group[\Bigg]{\myspace \sum_{z_k^j \in \Gamma_{\zeta}^r} (1-\abs{z_k^j}^2)^{1+\alpha}  }
    \end{align*} 
    as the number of elements of $z\0^j$ which lie inside $D(0,1/2)$ is bounded above uniformly in $j$ via \thref{?F?}. We now have
    $$
        \prod_{ 1/2 < \abs{z_k^j} < r_j } \frac{1}{\abs{z_k^j}^p} = \exp\group{ -p\log(\epsilon_j) D(z\0^j,0,r_j) } = \epsilon_j^{-p D(z\0^j,0,r_j) }
    $$
    by definition. Taking $\nu_j = \nu^{z\0^j}$, we get
    \begin{align*}
        (L/C)^p & \leq 
        \norm{f^j}_{T_{\alpha}^{p,q}(\nu_j)}^p \\&\leq
        C\epsilon_j^{-p D(z\0^j,0,r_j) } \int_{\bd \dd} \group[\Bigg]{ \myspace \sum_{z_k^j \in \Gamma_{\zeta}^r} (1-\abs{z_k^j}^2)^{1+\alpha}  }^{p/q} \dif l(\zeta) \\&\leq
        C'\epsilon_j^{p( (1+\alpha)/q - D(z\0^j,0,r_j) ) } \\&\leq
        C''\epsilon_j^{p( (1+\alpha)/q - D^-(z) ) } \epsilon_j^{-\epsilon_j}
    \end{align*}
    by \thref{?50?} and the lower sampling inequality. As $\epsilon_j \arrow 0$ and $\epsilon_j^{-\epsilon_j} \searrow 1$ for $\epsilon_j \leq 1/4$, we must have $p((1+\alpha)/q-D^-(z\0) ) \leq 0$. Thus we have $D^-(z\0) \geq (1+\alpha)/q$. 
    
    We can now choose some sequence $r_n \arrow 1$ and $w_n \in \dd$ such that $D(z\0\SP,w_n, r_n)$ decreases to $D^-(z\0)$. More concretely, we take 
    $$
        D^-(z\0)-1/n \leq D(z\0\SP,w_n,r_n) < D^-(z\0) + 1/n\,.
    $$
    Now set $z\0^n = g_{w_n}^{\eta}(z\0)$ for some $\eta>0$ sufficiently small and $g_{\dott}^{\eta}$ as in \thref{?P?}. As $\eta$ is sufficiently small, we have that $z\0^n$ is uniformly discrete and $z\0^n$ is still sampling by \thref{?P?,?48?}. Moreover, \thref{?P?} says that 
    \begin{align*}
        D^-(z\0) & \geq 
        D(z\0\SP, w_n, r_n) - 1/n \\&\geq
        (1-\eta)^{-1} \group{ D(z\0^n\,, w_n, r_n) - C \bracket{\log(1-r_n)}^{-1} } - 1/n \\& \geq
        (1-\eta)^{-1} \group{ D^-(z\0^n) - C \bracket{\log(1-r_n)}^{-1} } - 1/n
    \end{align*}
    for some $C$ depending only on $\delta$. Now, we apply the inequality above to get
    $$
        D^-(z\0^n) \geq (1+\alpha)/q
    $$
    which implies that
    $$
        D^-(z\0) \geq (1-\eta)^{-1} (1+\alpha)/q
    $$
    by letting $n$ go to infinity. Noting that $(1-\eta)^{-1} > 1$, we get the result.
\end{proof}

Combining \thref{?44?,?51?} gives \thref{samplingTheorem}, our second main theorem.

\newpage
\bibliography{ref.bib}

\end{document}